\documentclass[12pt]{article}
\usepackage{float}       
\usepackage{placeins}    
\usepackage{caption}     
\usepackage{tikz}
\usepackage{enumerate}
\usepackage[all,cmtip]{xy}
\usetikzlibrary{arrows,matrix,positioning}
\usepackage{soul}
\tikzset{>=stealth',
  head/.style = {fill = white, text=black},
  plaque/.style = {draw, rectangle, minimum size = 10mm},
  pil/.style={->,thick},
  junct/.style = {draw,circle,inner sep=0.5pt,outer sep=0pt, fill=black}
  }

\usepackage[colorlinks=true, pdfstartview=FitV, linkcolor=blue, citecolor=blue, urlcolor=blue]{hyperref}
\usepackage[english]{babel}
\usepackage{ytableau}
\usepackage[a4paper,top=2cm,bottom=2cm,left=3cm,right=3cm,marginparwidth=1.75cm]{geometry}

\usepackage{amssymb}
\usepackage{amsmath}
\usepackage[title]{appendix}
\usepackage{ntheorem}
\usepackage{siunitx}
\PassOptionsToPackage{hyphens}{url}\usepackage{hyperref}
\usepackage{cleveref}
\usepackage[utf8]{inputenc}
\usepackage{csquotes}
\usepackage{booktabs}
\usepackage{longtable}
\usepackage{adjustbox}
\usepackage{array}
\usepackage{url}
\usepackage{titlesec}
\usepackage{authblk}
\usepackage{xcolor} 

\newtheorem{theorem}{Theorem}[section]
\newtheorem{lemma}[theorem]{Lemma}

\newenvironment{myproof}
{\noindent\textit{Proof.}}
{\hfill$\blacksquare$}
\newtheorem{proposition}[theorem]{Proposition}
\newtheorem{corollary}[theorem]{Corollary}
\titleformat{\subsection}
  {\mdseries\itshape\large} 
  {\thesubsection}{1em}{} 

\crefformat{figure}{#2Figure~#1#3}
\Crefformat{figure}{#2Figure~#1#3}
\crefformat{table}{#2Table~#1#3}
\Crefformat{table}{#2Table~#1#3}
\crefformat{section}{#2Section~#1#3}
\Crefformat{section}{#2Section~#1#3}

\usepackage{soul}
\usepackage{xcolor}
\newtheorem{definition}[theorem]{Definition}
\newtheorem{remark}[theorem]{Remark}%
\newtheorem{problem}[theorem]{Problem}
\newtheorem{conjecture}[theorem]{Conjecture}
\newcommand{\ml}[1]{\mathbf{\color{red} #1} }
\newcommand{\bl}[1]{\mathbf{\color{blue} #1} }
\begin{document}

\def\l{{\lambda}}
\def\m{{\mu}}
\def\lm{{\l/\m}}
\def\qed{{\hfill$\blacksquare$}}
\newcommand{\bplus}{{\color{blue}+}}
\newcommand{\bminus}{{\color{red}-}}
\newcommand{\ZZ}{\mathbb{Z}}
\newcommand{\BB}{\mathbb{B}}
\newcommand{\NN}{\mathbb{N}}
\newcommand{\QQ}{\mathbb{Q}}
\newcommand{\RR}{\mathbb{R}}
\newcommand{\CC}{\mathbb{C}}
\newcommand{\gl}{\mathfrak{gl}}
\newcommand{\HW}{\mathrm{HW}}
\newcommand{\sr}{\mathrm{SVRPP}}
\newcommand{\ceq}{\mathrm{ceq}}
\newcommand{\ircont}{\mathrm{ircont}}
\newcommand{\ex}{\mathrm{ex}}
\newcommand{\flip}{\mathrm{flip}}
\newcommand{\res}{\mathrm{res}}
\newcommand{\bT} {\mathbb{T}}
\newcommand{\bS} {\mathbb{S}}
\newcommand{\bR} {\mathbb{R}}

\begin{center}
{\Large\bf Hybrid Grothendieck polynomials}

\vskip 6mm
{\small   }
Peter L. Guo, Mingyang Kang, Jiaji Liu

\end{center}

\begin{abstract}
For a skew shape $\lambda/\mu$, we define the hybrid Grothendieck polynomial 
$${G}_{\lambda/\mu}(\textbf{x};\textbf{t};\textbf{w})
=\sum_{T\in \mathrm{SVRPP}(\lambda/\mu)}
\textbf{x}^{\mathrm{ircont}(T)}\textbf{t}^{\mathrm{ceq} (T)}\textbf{w}^{\mathrm{ex}(T)}$$
as a weight generating function   over set-valued reverse plane partitions of shape $\lambda/\mu$. It   specializes to
\begin{itemize}
    \item[(1)] the  refined  
 stable Grothendieck polynomial introduced  by  Chan--Pflueger by setting  all $t_i=0$;

 \item[(2)] the refined  dual stable Grothendieck polynomial introduced  by Galashin--Grinberg--Liu by setting all $w_i=0$.
\end{itemize}

We show   that ${G}_{\lambda/\mu}(\textbf{x};\textbf{t};\textbf{w})$ is   symmetric    in the $\textbf{x}$ variables. By building  a crystal structure on set-valued reverse plane partitions, we obtain the  expansion of ${G}_{\lambda/\mu}(\textbf{x};\textbf{t};\textbf{w})$ in the basis of Schur functions, extending previous work by Monical--Pechenik--Scrimshaw and Galashin. Based on the Schur expansion, we deduce that hybrid Grothendieck polynomials of straight shapes have   saturated Newton polytopes. Finally, using Fomin--Greene's theory on  noncommutative  Schur functions, we give a combinatorial formula for  the image of ${G}_{\lambda/\mu}(\textbf{x};\textbf{t};\textbf{w})$ (in the case $t_i=\alpha$ and $w_i=\beta$) under the omega involution on symmetric functions. The formula unifies the structures of  weak set-valued tableaux and valued-set tableaux introduced by Lam--Pylyavskyy. 
Several   problems and conjectures are motivated and  discussed. 
\end{abstract}

\section{Introduction}
\def\e{\mathbf{e}}
\label{sec:introduction}

The stable Grothendieck polynomials represent the classes of    structure sheaves of Schubert varieties  in the K-theory of   Grassmannians. As discovered  by Buch  \cite{Buch}, stable Grothendieck polynomials admit a combinatorial model via set-valued   tableaux. As the Hopf-dual to stable Grothendieck polynomials with respect to the Hall inner product, the dual stable Grothendieck polynomials  represent the classes in the K-homology of  ideal sheaves
 of the boundaries of Schubert varieties, which  were described explicitly  by Lam and 
Pylyavskyy \cite{Lam}   using  reverse plane partitions. Expanding a dual stable Grothendieck polynomial of skew shape into   dual stable Grothendieck polynomials of straight shapes leads to a dual way to compute  the structure constants in the K-theory of   Grassmannians \cite{Lam,LMS}.
Refined versions of stable and dual stable Grothendieck polynomials have been explored by Chan  and Pflueger \cite{Chan} and Galashin, Grinberg and Liu \cite{Galashin 1}, respectively. We refer to  Section \ref{Last-22} for  various aspects of stable or dual stable Grothendieck polynomials. 

This paper provides  a simultaneous generalization of  stable and dual stable Grothendieck polynomials through introducing  {\it set-valued reverse plane partitions}. We   identify a partition $\lambda=(\lambda_1\geq \lambda_2\geq \cdots)$ with its  Young diagram, a left-justified array   with $\lambda_i$ boxes in row $i$. For partitions  $\mu\subseteq \lambda$, the skew diagram $\lambda/\mu$ is obtained from  $\lambda$ by removing the boxes in $\mu$.  For two nonempty  finite sets $A$ and $B$ of positive integers, write $A\leq B$ if $\max(A)\leq \min(B)$, and 
$A<B$ if $\max(A)< \min(B)$. 

\begin{definition}
A {\it set-valued reverse plane partition} (or abbreviated as SVRPP) of  shape $\lambda/\mu$ is a filling of nonempty finite sets of positive integers into the boxes of  $\lambda/\mu$  such that the sets are weakly increasing along each row and  each column.  
\end{definition}

\begin{figure}[h t]
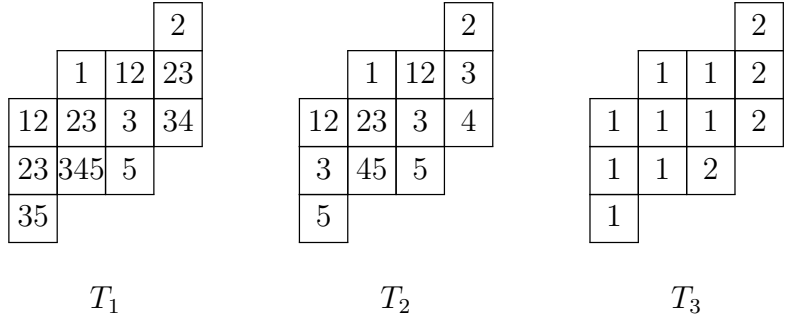

\begin{center}
\begin{tabular}{c c c c c c c}
 & & & & & & \\
\begin{ytableau}
\none&\none&\none&2\\
\none&1 & 12 &23\\
12  & 23 & 3 &34\\
23 & 345 & 5\\
35
\end{ytableau}& & &
\begin{ytableau}
\none&\none&\none&2\\
\none&1 & 12 & 3\\
12  & 23 & 3 & 4\\
3  & 45 & 5\\
5
\end{ytableau}& & &
\begin{ytableau}
\none&\none&\none&2\\
\none&1 & 1 &2\\
1  & 1 &  1 &2\\
1  & 1 &  2\\
1
\end{ytableau}\\
 & & & & & &\\
$T_1$ & & &$T_2$ & & &$T_3$ \\
\end{tabular}\\
\end{center}
\caption{$T_1$, $T_2$ and $T_3$ are  SVRPP's. Furthermore, $T_2$ is a set-valued tableau, and $T_3$ is a reverse plane  
partition. Here we have omitted  the braces $\{\ \}$ for sets as well as  the commas between numbers.}
\label{fig.1}
\end{figure}

Notice that a SVRPP recovers a {set-valued tableau} defined by Buch \cite{Buch} if the sets are strictly increasing in each column, and a {reverse plane partition} if each set contains   a single number. See Figure \ref{fig.1} for an illustration. 

We use $\mathrm{SVRPP}(\lambda/\mu)$ to denote the set of all SVRPP's of
shape $\lambda/\mu$.  Write $(i,j)$ for the box in row $i$ and column $j$ in the matrix coordinate, and for $T\in \mathrm{SVRPP}(\lambda/\mu)$, let $T(i,j)$ be the set filled in $(i,j)$.  Define 
\begin{itemize}
\item {\it irredundant
content}  
$\mathrm{ircont}(T) = (r_1,r_2, \ldots)$:   $r_i$
is the number of columns of $T$ that contain a set that has $i$;

\item  {\it column equalities vector} $\mathrm{ceq} (T) =
(c_1, c_2, \ldots)$: $c_i$ is the number of {\it redundant} boxes in $i$-th row  of $T$. Here a redundant box means a box $(i,j)$ such that $(i+1,j)\in \lambda/\mu$, and $\max{T(i,j)}=\min{T(i+1,j)}$;

\item {\it excess vector} $\mathrm{ex}(T)=(e_1, e_2, . . .)$: $e_i$ is equal to the sum of cardinalities  of subsets filled in the $i$-th row of $T$  minus  the number of  boxes
in the $i$-th row of $T$. That is, \[e_i=\sum_{(i,j)\in \lambda/\mu}\left(|T(i,j)|-1\right).\]
\end{itemize}
For the fillings in Figure \ref{fig.1}, we have 
\begin{align*}
\mathrm{ircont}(T_1)=(3,4,4,2,3), \ \ \  \mathrm{ceq} (T_1)=(1,1,2,1),   \ \ \  \mathrm{ex} (T_1)=(0,2,3,3,1). \\[5pt]
\mathrm{ircont}(T_2)=(3,4,4,2,3), \ \ \  \mathrm{ceq} (T_2)=(0,0,0,0),   \ \ \  \mathrm{ex} (T_2)=(0,1,2,1,0). \\[5pt]
\mathrm{ircont}(T_3)=(3,2,0,0,0), \ \ \  \mathrm{ceq} (T_3)=(1,3,2,1),   \ \ \  \mathrm{ex} (T_3)=(0,0,0,0,0).
\end{align*}

We define the {\it hybrid Grothendieck polynomial} of shape $\lambda/\mu$ as   
\begin{equation}\label{main-func}
  {G}_{\lambda/\mu}(\textbf{x};\textbf{t};\textbf{w})
=\sum_{T\in \mathrm{SVRPP}(\lambda/\mu)}
\textbf{x}^{\mathrm{ircont}(T)}\textbf{t}^{\mathrm{ceq} (T)}\textbf{w}^{\mathrm{ex}(T)},  
\end{equation}
where $\textbf{x}=(x_1,x_2,\ldots)$, 
$\textbf{t}=(t_1,t_2,\ldots)$ and $\textbf{w}=(w_1,w_2,\ldots)$ are coutable sequences of variables. Here, for a weak composition $\alpha=(\alpha_1,\alpha_2,\ldots)\in \mathbb{Z}_{\geq 0}^n$ and a sequence $\mathbf{y}=(y_1,y_2,\ldots)$ of variables, we adopt the convention $\mathbf{y}^\alpha=y_1^{\alpha_1}y_2^{\alpha_2}\cdots$.

Our first  result  reveals  that   ${G}_{\lambda/\mu}({\mathbf{x}} ;\mathbf{t} ; \mathbf{w})$ is symmetric in the $\mathbf{x}$ variables.

\begin{theorem}\label{THM-1}
For any skew shape $\lambda/\mu$, ${G}_{\lambda/\mu}({\mathbf{x}} ;\mathbf{t} ; \mathbf{w})$ is symmetric in $\mathbf{x}$.    
\end{theorem}

Hybrid Grothendieck polynomials unify stable and dual stable Grothendieck polynomials.

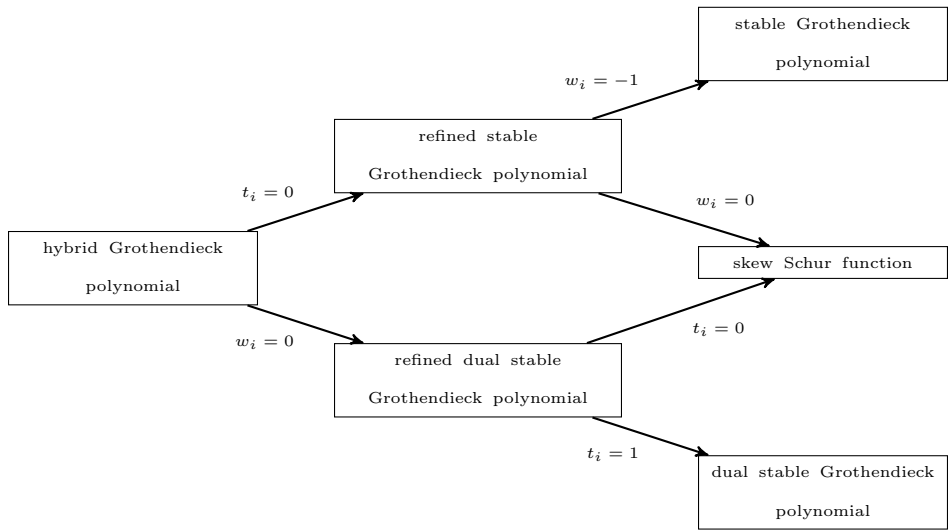
\begin{figure}[h t]
    \centering
    \begin{tikzpicture}
  \node(A)[draw, rectangle, text width=3cm, align=center] {\tiny hybrid Grothendieck polynomial};
  \node[above right = 0.5 and 1 of A](B)[draw, rectangle, text width=3.5cm, align=center] {\tiny refined stable Grothendieck polynomial};
  \node[below right = 0.5 and 1 of A](C)[draw, rectangle, text width=3.5cm, align=center] {\tiny refined dual stable Grothendieck polynomial};

  \node[above right = 0.5 and 1 of B](D)[draw, rectangle, text width=3cm, align=center] { \tiny stable Grothendieck polynomial};

   \node[below right = 0.7 and 1 of B](E)[draw, rectangle, text width=3cm, align=center] {\tiny    skew  Schur function   };

   \node[below right = 0.5 and 1 of C](F)[draw, rectangle, text width=3cm, align=center] {{\tiny dual stable Grothendieck polynomial}};

   \path (A) edge[pil,color=black] node[above left,black]{\tiny $t_i=0$} (B);

   \path (A) edge[pil,color=black] node[below left,black]{\tiny $w_i=0$} (C);

   \path (B) edge[pil,color=black] node[above left,black]{\tiny $w_i=-1$} (D);

   \path (B) edge[pil,color=black] node[above right,black]{\tiny $w_i=0$} (E);

   \path (C) edge[pil,color=black] node[below right,black]{\tiny $t_i=0$} (E);

   \path (C) edge[pil,color=black] node[below left,black]{\tiny $t_i=1$} (F);
\end{tikzpicture}
    \caption{Specializations of hybrid Grothendieck polynomials.}
    \label{fig:my_label-yt}
\end{figure}

\begin{itemize}

\item   $t_i=0$:  In this case, only set-valued tableaux have contributions, and thus  ${G}_{\lambda/\mu}(\mathbf{x};\mathbf{0};\mathbf{w})$ is   the refined stable Grothendieck polynomial  defined  by Chan  and Pflueger \cite{Chan}. If we further set all $w_i=-1$, then ${G}_{\lambda/\mu}(\textbf{x};\textbf{0};\mathbf{-1})$ becomes  the ordinary stable Grothendieck polynomial  \cite{Buch}.

We point out  that the precise form of a refined stable Grothendieck polynomial in \cite{Chan} is
${G}_{\lambda/\mu}(\mathbf{x};\mathbf{0};-\mathbf{w})$, where $-\mathbf{w}=(-w_1, -w_2,\ldots)$. For notational simplicity, we adopt the current description as above. 

\item   $w_i=0$:  In this case, only ordinary reverse plane partitions have contributions, and thus  ${G}_{\lambda/\mu}(\textbf{x};\textbf{t};\textbf{0})$ is  the refined dual stable Grothendieck polynomial  introduced  by Galashin, Grinberg and Liu \cite{Galashin 1}. If we further set all $t_i=1$, then ${G}_{\lambda/\mu}(\textbf{x};\textbf{1};\textbf{0})$ is   the ordinary dual stable Grothendieck polynomial due to Lam and 
Pylyavskyy \cite{Lam}.

\end{itemize}
Of course, when both $t_i$ and $w_i$ are equal to zero, only semistandard Young tableaux have contributions, and hence  ${G}_{\lambda/\mu}({\mathbf{x}} ;\mathbf{0} ; \mathbf{0})$ specializes to the skew Schur function $s_{\lambda/\mu}(\mathbf{x})$. The above specializations are depicted in Figure  \ref{fig:my_label-yt}.

The proof of Theorem \ref{THM-1} is in the spirit of the method in  \cite{Galashin 1}. To do this,  we show that for  any positive integer $i$, there  exists
an involution $\Phi$  on $\mathrm{SVRPP} (\lambda/\mu)$ which preserves both the ceq statistic  and the ex statistic, 
while switching  the  $i$-th entry and $(i + 1)$-th
entry in the ircont statistic. When restricted to ordinary reverse plane partitions, $\Phi$ coincides with the involution constructed   in \cite{Galashin 1}. 

Our second   result is 
to expand  a hybrid Grothendieck polynomial in terms of Schur polynomials. Let $\mathrm{SVRPP}^n(\lm)$ denote the collection of SVRPP's of shape $\lambda/\mu$ filled   with nonempty subsets of  $[n]=\{1,2,\ldots,n\}$. Based  on the algorithm  in  the construction of $\Phi$, we access  a $\gl_n$-crystal structure on $\mathrm{SVRPP}^n(\lm)$. We  show  that each connected component in its crytal graph is isomorphic to the crystal  $\mathcal{B}(\nu)$ (for some partition $\nu$) corresponding to the crystal basis of the highest module $V (\nu)$ of the quantum group $U_q(\gl_n)$. As a consequence, we obtain the following  Schur expansion.

\begin{theorem}\label{THM-2}
Let   $\mathbf{x}_n=(x_1,\ldots, x_n)$. 
 For any skew shape $\lambda/\mu$, we have \begin{equation}\label{expan}
    {G}_{\lm}(\mathbf{x}_n ;\mathbf{t};\mathbf{w})=\sum_{\gamma}\sum_{\theta}\mathbf{t}^\gamma\mathbf{w}^\theta\sum_{\nu}H_{\lm,n}^{\nu,\gamma,\theta}s_{\nu}(\mathbf{x}_n),
\end{equation}
 where $\gamma$ and $\theta$ are  over weak compositions, $\nu $ is over   partitions, and  the coefficient $H_{\lm,n}^{\nu,\gamma,\theta}$ counts the   number of  $T\in \mathrm{SVRPP}^n(\lambda/\mu)$   with  $\mathrm{ircont}(T)=\nu$, $\mathrm{ceq}(T) = \gamma$ and $\mathrm{ex}(T)=\theta$,  such that the reading word of $T$  
is a  reverse lattice word.
\end{theorem}

Theorem \ref{THM-2} tells that hybrid Grothendieck polynomials are Schur positive. 
Let us take $\lambda=(2,2)$, $\mu=(1)$ and $n=3$. 
Appendix \ref{hgi-6-9} lists all the connected components in the crystal graph of  $\mathrm{SVRPP}^3(\lambda/\mu)$. By direct computation, we have 
\begin{align}
{G}_{\lm}(\mathbf{x}_3;\mathbf{t};\mathbf{w})&=
t_1s_{(2)}(\mathbf{x}_3)+\left(1+t_1w_2\right)s_{(2,1)}(\mathbf{x}_3)+\left(w_2+t_1w_1+t_1w_1w_2\right)s_{(2,2)}(\mathbf{x}_3)\nonumber\\[5pt]
&\ \ + \left(w_1+w_2+t_1w_1^2+t_1w_1w_2+t_1w_2^2\right)s_{(2,1,1)}(\mathbf{x}_3)\nonumber\\[5pt]
&\ \ +\left(w_1w_2+w_2^2+t_1w_1^2w_2+t_1w_1w_2^2\right)s_{(2,2,1)}(\mathbf{x}_3)\nonumber\\[5pt]
&\ \ +\left(w_1w_2^2+t_1w_1^2w_2^2\right)s_{(2,2,2)}(\mathbf{x}_3).\label{Iy-0-6}
\end{align}

We point out some  specializations of Theorem \ref{THM-2}. 
\begin{itemize}
    \item  $t_i=0$: In this case,   \eqref{expan} reduces to  the expansion of the  refined stable Grothendieck polynomial ${G}_{\lambda/\mu}(\textbf{x}_n;\textbf{0};\textbf{w})$.   If  further setting all $w_i=-1$ and $\mu=\emptyset$, then we get the expansion of  a stable Grothendieck polynomial   of straight shape. The expansion for straight shape was first  derived   by Lenart \cite{Lenart} in an algebraic manner. Buch  \cite{Buch} gave  a combinatorial proof by introducing the uncrowding  algorithm. In fact, Lenart's formula was expressed in terms of   increasing tableaux, which were shown  by  Monical,  Pechenik,  Scrimshaw   \cite{Monical} to be equivalent to the language of  lattice words.

\item $w_i=0$: In this case, \eqref{expan} becomes the expansion  of the refined dual Grothedieck polynomial   ${G}_{\lambda/\mu}(\textbf{x}_n;\textbf{t};\textbf{0})$  into Schur polynomials by Galashin \cite{Galashin 2}. 
A tabloid crystal  perspective for the Schur expansion of 
${G}_{\lambda/\mu}(\textbf{x}_n;\textbf{1};\textbf{0})$ appeared in the work of Li,  Morse and   Shields \cite{LMS}. 
The Schur expansion of ${G}_{\lambda}(\textbf{x}_n;\textbf{1};\textbf{0})$  of straight shape was first given by Lam and 
Pylyavskyy \cite{Lam}.
    
  \item   $t_i=w_i=0$: In this particular case,  \eqref{expan} is the expansion of the skew Schur polynomial  $s_{\lambda/\mu}(\mathbf{x}_n)$, which recovers  the  classic Littlewood--Richardson coefficients, namely, the structure constants in the cohomology ring  of   Grassmannians . 
\end{itemize}

We apply   Theorem \ref{THM-2} to study  the Newton polyotpe  of  ${G}_{\lambda}(\mathbf{x}_n ;\mathbf{t};\mathbf{w})$ of straight shape.  Recall that the Newton polytope $\mathrm{Newton}(f(\mathbf{x}_n))$ of a polynomial \[f(\mathbf{x}_n)=\sum_{\gamma\in \mathbb{Z}_{\geq 0}^n}c_\gamma\, \mathbf{x}_n^\gamma\]
is the convex hull of the set $\{\gamma\colon c_\gamma \neq 0\}$ of its exponent vectors. Clearly,   each exponent vector appears as a lattice point in $\mathrm{Newton}(f(\mathbf{x}_n))$. If the converse is true, that is, each lattice point in $\mathrm{Newton}(f(\mathbf{x}_n))$ is an exponent vector of $f(\mathbf{x}_n)$, then we say that $f(\mathbf{x}_n)$ has  saturated Newton polytope (SNP). 
Since the pioneer work of Monical, Tokcan and  Yong \cite{Yong-1}, the SNP phenomenon in algebraic combinatorics  has been widely investigated, see for example \cite{Besson,Castillo-1,Castillo-2,Esbobar,Fink,Huh,Matherne,KSS,Nguyen}.

Escobar and Yong \cite{Esbobar}   showed   that any stable Grothendieck polynomial   of straight shape has SNP. Nguyen,  Ngoc,   Tuan and Do Le Hai   \cite{Nguyen} confirmed the SNP property of any dual stable Grothendieck polynomial  of straight shape.

\begin{theorem}\label{Newton-98}
For any  partition $\lambda$,   $G_{\l}(\mathbf{x}_n;\mathbf{t};\mathbf{w})$  has SNP. Here $\mathbf{t}$ and $\mathbf{w}$ are considered as generic parameters, or in other words, we may assume all $t_i=w_i=1$.
\end{theorem}

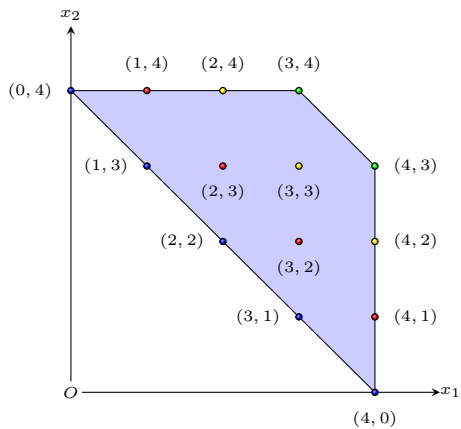
\begin{figure}[h t]    
\begin{center}
\begin{tikzpicture}[scale=1]
	\def\a{blue}
	\def\b{red}
	\def\c{pink}
	\def\d{violet}
	\def\e{orange}
	\def\opacity{100}
	\tikzstyle{point1}=[ball color=blue, circle, draw=black, inner sep=0.03cm]
	\tikzstyle{point2}=[ball color=red, circle, draw=black, inner sep=0.03cm]
	\tikzstyle{point3}=[ball color=yellow, circle, draw=black, inner sep=0.03cm]
    \tikzstyle{point4}=[ball color=green, circle, draw=black, inner sep=0.03cm]

    \node (O0) at (0,0)[]{};
	\node (O1) at (5,0)[]{};
	\node (O2) at (0,5)[]{};
    
    \node at (O0) []{\tiny$O$}; 
	\node at (O1) []{\tiny$x_1$}; 
	\node at (O2) []{\tiny$x_2$};
 
    \draw[-{stealth[scale=3.0]}] (O0) -- (O1);
	\draw[-{stealth[scale=3.0]}] (O0) -- (O2);

\filldraw[fill=\a!20,rounded corners=0.5pt] (0,4) -- (1,4) -- (2,4) -- (3,4) -- (4,3) -- (4,2) -- (4,1) -- (4,0) -- (3,1) -- (2,2) -- 
(1,3)--cycle;

    \node (A1) at (0,4)[point1]{};
	\node (A2) at (1,3)[point1]{};
	\node (A3) at (2,2)[point1]{};
	\node (A4) at (3,1)[point1]{};
	\node (A5) at (4,0)[point1]{};

    \node at (A1) [left = 1mm]{\tiny$(0,4)$}; 
	\node at (A2) [left = 1mm]{\tiny$(1,3)$};
	\node at (A3) [left = 1mm]{\tiny$(2,2)$};
	\node at (A4) [left = 1mm]{\tiny$(3,1)$};
	\node at (A5) [below = 1mm]{\tiny$(4,0)$};

    \node (B1) at (1,4)[point2]{};
	\node (B2) at (2,3)[point2]{};
	\node (B3) at (3,2)[point2]{};
    \node (B4) at (4,1)[point2]{};

    \node at (B1) [above = 1mm]{\tiny$(1,4)$}; 
	\node at (B2) [below = 1mm]{\tiny$(2,3)$};
	\node at (B3) [below = 1mm]{\tiny$(3,2)$};
    \node at (B4) [right=1mm]{\tiny$(4,1)$};

    \node (C1) at (2,4)[point3]{};
	\node (C2) at (3,3)[point3]{};
	\node (C3) at (4,2)[point3]{};

    \node at (C1) [above = 1mm]{\tiny$(2,4)$}; 
	\node at (C2) [below = 1mm]{\tiny$(3,3)$};
	\node at (C3) [right = 1mm]{\tiny$(4,2)$};

    \node (D1) at (3,4)[point4]{};
	\node (D2) at (4,3)[point4]{};

    \node at (D1) [above = 1mm]{\tiny$(3,4)$}; 
	\node at (D2) [right = 1mm]{\tiny$(4,3)$};
\end{tikzpicture}
\vspace{-15pt}
\end{center}
\caption{The Newton polytope  of $G_{(4,2,1)}(\mathbf{x}_2;\mathbf{1};\mathbf{1})$. }
    \label{fig:my_label=7-3}
\end{figure}

For example, we have  
\begin{align*}
G_{(4,2,1)}(\mathbf{x}_2;\mathbf{1};\mathbf{1})
&=s_{(4)}(\mathbf{x}_2)+5s_{(4,1)}(\mathbf{x}_2)+10s_{(4,2)}(\mathbf{x}_2)+6s_{(4,3)}(\mathbf{x}_2)\\[5pt]
&=\left(x^{(4,0)}+x^{(3,1)}+x^{(2,2)}+x^{(1,3)}+x^{(0,4)}\right)+5\left(x^{(4,1)}+x^{(3,2)}+x^{(2,3)}\right.\\
&\ \ \ \ \left.+x^{(1,4)}\right)+10\left(x^{(4,2)}+x^{(3,3)}+x^{(2,4)}\right)+6\left(x^{(4,3)}+x^{(3,4)}\right),
\end{align*}
which has Newton polytope as drawn in Figure \ref{fig:my_label=7-3}.

To prove  Theorem  \ref{Newton-98}, we apply    Theorem \ref{THM-2} to show that for any given $\lambda$, there   exist partitions $\mu\subseteq \nu$ such that $s_\rho(\mathbf{x}_n)$ appears in the expansion of $G_{\l}(\mathbf{x}_n;\mathbf{t};\mathbf{w})$ if and only if  $\mu \subseteq \rho\subseteq \nu$.
This, along with  a criterion for SNP given in \cite{Nguyen}, leads to a proof of Theorem  \ref{Newton-98}.  As a byproduct, we get a combinatorial description for the degree of $G_{\l}(\mathbf{x}_n;\mathbf{t};\mathbf{w})$, see Corollary \ref{bnf-5-8}.

Finally, we consider the image of a hybrid Grothendieck polynomial under the {\it omega involution} $\omega$ on symmetric functions. Recall that $\omega$ can be defined by sending a Schur function $s_\lambda(\textbf{x})$ to $s_{\lambda'}(\textbf{x})$, where $\lambda'$ is the conjugate of $\lambda$. In the case of setting all $t_i=\alpha$ and $w_i=\beta$, we employ   Fomin--Greene type operators \cite{Fomin} to generate 
\[
{G}_{\lm}(\mathbf{x} ;\alpha;\beta):={G}_{\lm}(\mathbf{x} ;\mathbf{t};\mathbf{w})|_{t_i=\alpha, w_i=\beta}.
\] 
This enables us to capture  the image, denote ${J}_{\lm}(\mathbf{x} ;\alpha;\beta)$, of ${G}_{\lm}(\mathbf{x} ;\alpha;\beta)$ under the involution $\omega$.  

\begin{theorem}\label{THM-33}
For any skew shape $\lambda/\mu$, we have 
\begin{equation}\label{Hi-0-3}
 {J}_{\lm}(\mathbf{x} ;\alpha;\beta)=\sum_{T}   \alpha^{|\mathrm{mark}(T)|}\beta^{|\mathrm{ex}(T)|} \mathbf{x}^{\mathrm{umcont}(T)},
\end{equation}
where the sum is over all 
 marked multiset-valued  tableaux (see Section \ref{info-14} for the precise definition).   
\end{theorem}

The   formula for  ${J}_{\lm}(\mathbf{x} ;\alpha;\beta)$  in \eqref{Hi-0-3} unifies  the images of stable and dual stable Grothendieck polynomials established  by   Lam and 
Pylyavskyy \cite{Lam}. Precisely, 
\begin{itemize}
    \item $\alpha=0$ and $\beta=-1$:  ${J}_{\lm}(\mathbf{x} ;0;-1)$ is equal to the image of the stable Grothendieck polynomial ${G}_{\lm}(\mathbf{x} ;\mathbf{0};-\mathbf{1})$;

    \item $\alpha=1$ and $\beta=0$:  ${J}_{\lm}(\mathbf{x} ;1;0)$ is equal to the image of the dual stable Grothendieck polynomial   ${G}_{\lm}(\mathbf{x} ;\mathbf{1};\mathbf{0})$. 
\end{itemize}


This paper is structured as follows. Sections \ref{SC-11}, \ref{Sc-43}, \ref{Sec-5--2} and \ref{info-14} are devoted to proofs of Theorems \ref{THM-1},  \ref{THM-2}, \ref{Newton-98} and \ref{THM-33}, respectively. In Section \ref{Last-22}, we conclude with  several    problems and conjectures that are motivated by our results. 

\subsection*{Acknowledgements}

We are grateful to  Rui Xiong  and Tianyi Yu for very helpful comments and suggestions. 
This work was  supported by the National Natural Science Foundation of China (No. 12371329) and the Fundamental Research Funds for the Central Universities (No. 63243072).

\section{Proof of Theorem \ref{THM-1}}\label{SC-11}

Theorem \ref{THM-1} is equivalent to  the following   statement.  

\begin{theorem}\label{UHI-1}
For any skew shape $\lambda/\mu$ and any positive integer  $i $, there exists
an involution $\Phi_i$ on $\mathrm{SVRPP} (\lambda/\mu) $    which preserves the $\mathrm{ceq}$ and $\mathrm{ex}$  statistics,
whereas acts on the $\mathrm{ircont}$   statistic by exchanging  its $i$-th and $(i + 1)$-th
entries.    
\end{theorem}

The construction of $\Phi_i$ can be regarded as a set-valued analogue of that given by Galashin, Grinberg and Liu \cite{Galashin 1}.  Indeed,  the approach   in \cite{Galashin 1} is our original motivation for searching for a proof of Theorem \ref{THM-1}.  Since the involution $\Phi_i$ concerns only the entries $i$ and $i+1$, we may pay attention to the set of  two-element fillings, that is, the   set $\mathrm{SVRPP}^{n} (\lambda/\mu)$ for $n=2$. 
In order to emphasize that each box is filled with a subset of $\{1,2\}$, we use  $\mathrm{SVRPP}^{12} (\lambda/\mu)$ instead of $\mathrm{SVRPP}^{2} (\lambda/\mu)$.




\begin{theorem}\label{HGT-22}
There exists an involution $\Phi $  on $\mathrm{SVRPP}^{12} (\lambda/\mu)$ which
preserves the $\mathrm{ceq}$ and $\mathrm{ex}$  statistics, whereas switches the first two entries in the $\mathrm{ircont}$ statistic.    
\end{theorem}

Once the involution $\Phi$ in Theorem  \ref{HGT-22} is established, we can define $\Phi_i(T)$ for $T\in \mathrm{SVRPP}(\lambda/\mu)$  by applying $\Phi$ to the subtabeau $T^{(i)}$ of $T$ occupied by  $i$ and $i+1$  (when doing $\Phi$, we treat $i$ and $i+1$ as 1 and 2, respectively).

\subsection{Construction of the involution  \texorpdfstring{$\Phi$}{Phi}}
\label{sec:method}

Let us begin by defining a  {\it 12-SV table} as a filling $T$ of $\lambda/\mu$ with nonempty subsets of $\{1,2\}$   such that the sets in each column are weakly increasing (the sets in each row are not necessarily weakly increasing). Clearly, each SVRPP is a 12-SV table. 
A  column in a 12-SV table is called {\it 1-pure} if it contains only $\{1\}$’s, and {\it 2-pure} if it contains only  $\{2\}$’s, and {\it mixed}  otherwise.

The {\it flip} operation on   12-SV tables is defined as follows. For a 12-SV table $T$, $\mathrm{flip}(T)$ is obtained  by changing each column of $T$ in the following way:  
(i) for a 1-pure column, replace each $\{1\}$ by $\{2\}$; 
(ii) for a 2-pure column, replace each $\{2\}$ by $\{1\}$;   
 (iii) keep each mixed column unchanged.

 A column index $k$ is called a \textit{descent} of a 12-SV table $T$ if there exist two  boxes  $(i,k ),
(i, k + 1) \in \lambda/\mu$ such that $\max( T (i, k) )= 2$ and $\min (T (i, k + 1)) = 1$. Evidently, a 12-SV table   is a 12-SVRPP if and only if it has  no descent.
According to the  configurations  of columns $k$ and $k + 1$,   descents can be classified into four categories:
\begin{itemize}
    \item  (M1): column   $k$ is mixed and column $(k + 1)$ is $1$-pure;

\item (2M): column   $k$ is 2-pure and column $(k + 1)$ is mixed;

\item (21): column   $k$ is 2-pure and column $(k + 1)$ is $1$-pure;

\item (MM):  both columns $k$ and $k+1$ are mixed.
\end{itemize}

\begin{figure}[h t]
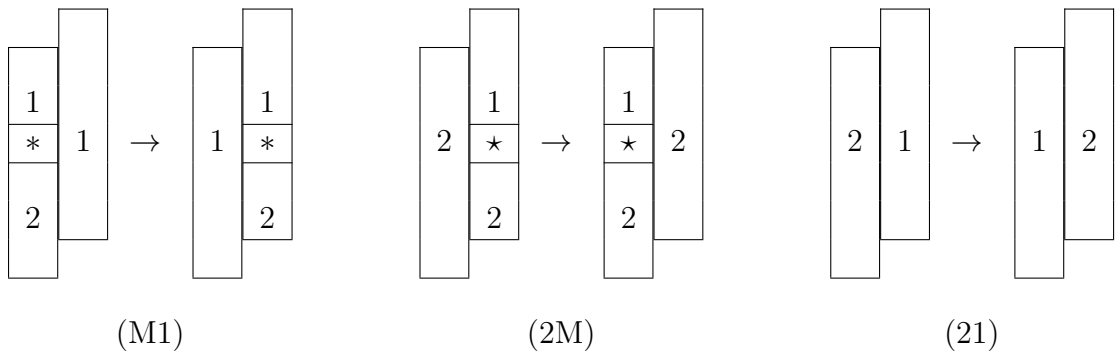


\begin{tabular}{c c c c c }
 &      &      &     & \\
\begin{tabular}{cc}
\cline{2-2} & \multicolumn{1}{|c|}{}\\
\cline{1-1} \multicolumn{1}{|c|}{} & \multicolumn{1}{|c|}{} \\
\multicolumn{1}{|c|}{1} & \multicolumn{1}{|c|}{}\\
\cline{1-1} \multicolumn{1}{|c|}{$*$} & \multicolumn{1}{|c|}{1}\\
\cline{1-1}\multicolumn{1}{|c|}{} & \multicolumn{1}{|c|}{}\\
\multicolumn{1}{|c|}{2} & \multicolumn{1}{|c|}{}\\
\cline{2-2} \multicolumn{1}{|c|}{} \\
\cline{1-1}
\end{tabular}

\ $\rightarrow$\quad

\begin{tabular}{cc}
\cline{2-2} & \multicolumn{1}{|c|}{}\\
\cline{1-1} \multicolumn{1}{|c|}{} & \multicolumn{1}{|c|}{} \\
\multicolumn{1}{|c|}{} & \multicolumn{1}{|c|}{1}\\
\cline{2-2} \multicolumn{1}{|c|}{1} & \multicolumn{1}{|c|}{$*$}\\
\cline{2-2}\multicolumn{1}{|c|}{} & \multicolumn{1}{|c|}{}\\
\multicolumn{1}{|c|}{} & \multicolumn{1}{|c|}{2}\\
\cline{2-2} \multicolumn{1}{|c|}{} \\
\cline{1-1}
\end{tabular}
&\quad\quad\quad&
\begin{tabular}{cc}
\cline{2-2} & \multicolumn{1}{|c|}{}\\
\cline{1-1} \multicolumn{1}{|c|}{} & \multicolumn{1}{|c|}{} \\
\multicolumn{1}{|c|}{} & \multicolumn{1}{|c|}{1}\\
\cline{2-2} \multicolumn{1}{|c|}{2} & \multicolumn{1}{|c|}{$\star$}\\
\cline{2-2} \multicolumn{1}{|c|}{} & \multicolumn{1}{|c|}{}\\
\multicolumn{1}{|c|}{} & \multicolumn{1}{|c|}{2}\\
\cline{2-2} \multicolumn{1}{|c|}{} \\
\cline{1-1}
\end{tabular}

\ $\rightarrow$\quad

\begin{tabular}{cc}
\cline{2-2} & \multicolumn{1}{|c|}{}\\
\cline{1-1} \multicolumn{1}{|c|}{} & \multicolumn{1}{|c|}{} \\
\multicolumn{1}{|c|}{1} & \multicolumn{1}{|c|}{}\\
\cline{1-1} \multicolumn{1}{|c|}{$\star$} & \multicolumn{1}{|c|}{2}\\
\cline{1-1} \multicolumn{1}{|c|}{} & \multicolumn{1}{|c|}{}\\
\multicolumn{1}{|c|}{2} & \multicolumn{1}{|c|}{}\\
\cline{2-2} \multicolumn{1}{|c|}{} \\
\cline{1-1}
\end{tabular}
&\quad\quad\quad&
\begin{tabular}{cc}
\cline{2-2} & \multicolumn{1}{|c|}{}\\
\cline{1-1} \multicolumn{1}{|c|}{} & \multicolumn{1}{|c|}{} \\
\multicolumn{1}{|c|}{} & \multicolumn{1}{|c|}{}\\
\multicolumn{1}{|c|}{2} & \multicolumn{1}{|c|}{1}\\
\multicolumn{1}{|c|}{} & \multicolumn{1}{|c|}{}\\
\multicolumn{1}{|c|}{} & \multicolumn{1}{|c|}{}\\
\cline{2-2} \multicolumn{1}{|c|}{} \\
\cline{1-1}
\end{tabular}

\ $\rightarrow$\quad

\begin{tabular}{cc}
\cline{2-2} & \multicolumn{1}{|c|}{}\\
\cline{1-1} \multicolumn{1}{|c|}{} & \multicolumn{1}{|c|}{} \\
\multicolumn{1}{|c|}{} & \multicolumn{1}{|c|}{}\\
\multicolumn{1}{|c|}{1} & \multicolumn{1}{|c|}{2}\\
\multicolumn{1}{|c|}{} & \multicolumn{1}{|c|}{}\\
\multicolumn{1}{|c|}{} & \multicolumn{1}{|c|}{}\\
\cline{2-2} \multicolumn{1}{|c|}{} \\
\cline{1-1}
\end{tabular}\\
&  &  &  &\\
 (M1) &  &  (2M)  &  &  (21)  \\
\end{tabular}

\caption{The descent-resolution operations. } \label{fig:des-res-preliminary}
\end{figure}

\begin{figure}[ht]
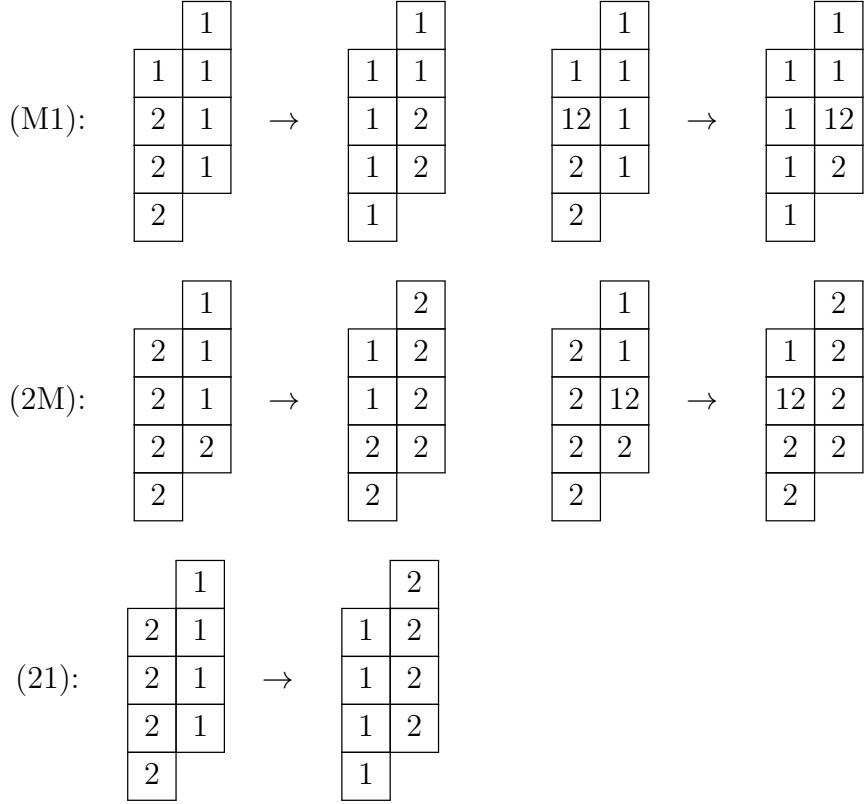

    \centering
    \begin{tabular}{cc}
    (M1):\quad
    \begin{tabular}{c}
    \begin{ytableau}
    \none & 1\\
    1& 1\\
    2& 1\\
    2& 1\\
    2
    \end{ytableau}     
    \end{tabular}
        
    \ $\rightarrow$\quad
    
    \begin{tabular}{c}
    \begin{ytableau}
    \none & 1\\
    1& 1\\
    1& 2\\
    1& 2\\
    1
    \end{ytableau}   
    \end{tabular}
    \quad\quad
    & 
    \begin{tabular}{c}
    \begin{ytableau}
    \none & 1\\
    1& 1\\
    12& 1\\
    2& 1\\
    2
    \end{ytableau}     
    \end{tabular}
        
    \ $\rightarrow$\quad
    
    \begin{tabular}{c}
    \begin{ytableau}
    \none & 1\\
    1& 1\\
    1& 12\\
    1& 2\\
    1
    \end{ytableau}   
    \end{tabular}
     \\
     & \\
     (2M):\quad
    \begin{tabular}{c}
    \begin{ytableau}
    \none & 1\\
    2& 1\\
    2& 1\\
    2& 2\\
    2
    \end{ytableau}     
    \end{tabular}
        
    \ $\rightarrow$\quad
    
    \begin{tabular}{c}
    \begin{ytableau}
    \none & 2\\
    1& 2\\
    1& 2\\
    2& 2\\
    2
    \end{ytableau}   
    \end{tabular}
     \quad\quad&
    \begin{tabular}{c}
    \begin{ytableau}
    \none & 1\\
    2& 1\\
    2& 12\\
    2& 2\\
    2
    \end{ytableau}     
    \end{tabular}
        
    \ $\rightarrow$\quad
    
    \begin{tabular}{c}
    \begin{ytableau}
    \none & 2\\
    1& 2\\
    12& 2\\
    2& 2\\
    2
    \end{ytableau}   
    \end{tabular}\\
    &\\
     (21):\quad
    \begin{tabular}{c}
    \begin{ytableau}
    \none & 1\\
    2& 1\\
    2& 1\\
    2& 1\\
    2
    \end{ytableau}     
    \end{tabular}
     \ $\rightarrow$\quad
    
    \begin{tabular}{c}
    \begin{ytableau}
    \none & 2\\
    1& 2\\
    1& 2\\
    1& 2\\
    1
    \end{ytableau}   
    \end{tabular}
    \quad\quad&
    
    \end{tabular}
    \caption{Examples of the descent-resolution operations.}
    \label{fig:my_label-2-78}
\end{figure}

As will been seen later, we shall {\it not} encounter the (MM) situation in the construction of $\Phi$. 
For each of the first three kinds of descents, we will    define a descent-resolution operation, 
denoted  $ \mathrm{res}_k$, to resolve
the descent by changing
the entries in the $k$-th and ($k + 1$)-th columns. The three kinds of descent-resolution operations are  illustrated in   Figure \ref{fig:des-res-preliminary}, where  the two  $*$'s in (M1) are either both equal to $\{2\}$ or both equal to $\{1,2\}$, while   the   two $\star$'s in (2M) are either both equal to $\{1\}$ or both equal to $\{1,2\}$. See Figure  \ref{fig:my_label-2-78} for some concrete examples. 

\begin{remark}
When neither   $*$ nor $\star$ is equal to $\{1,2\}$,    Figure  \ref{fig:des-res-preliminary} reduces to  the descent-resolution  operations  given by  Galashin, Grinberg and Liu \cite{Galashin 1}. 
\end{remark}

The map $\Phi$ on $\mathrm{SVRPP}^{12} (\lambda/\mu)$ can now be described as follows.
We   use   
12-SVRPP to represent a SVRPP in $\mathrm{SVRPP}^{12} (\lambda/\mu)$.  
\begin{itemize}
    \item Let $T$ be a 12-SVRPP. First, apply the flip map  to $T$  to
get a 12-SV table $\mathrm{flip}(T)$. Then, apply the descent-resolution operations  in  Figure \ref{fig:des-res-preliminary} to $\mathrm{flip}(T)$ (in arbitrary
order) until we arrive at a 12-SV table with no descents, say $P$. Set $\Phi(T)= P$.
\end{itemize}

In the remaining subsections, we will show  that  the map $\Phi$ is well defined.  Roughly speaking, we need to check that  the process
terminates after a finite steps of descent-resolution operations, and   $\Phi(T)$ is independent of the order of applying the operations. Moreover, we have to  verify  that $\Phi$ is an involution  satisfying  the requirements in   Theorem \ref{HGT-22}.

In Figure \ref{VBIO-1}, we give an illustration of the construction of $\Phi$. We
start with a 12-SVRPP $T$, and then apply the flip map followed by a sequence of descent-resolution operations. Eventually, we get a  12-SVRPP which is defined as    $\Phi(T)$. 
\begin{figure}[ht]
\begin{tabular}{ccc}
& &  \\
\begin{ytableau}
\none&\none&\none&2\\
\none&1 & 1 &2\\
1  & 1 & 12 &2\\
1  & 1 & 2  &2\\
1  & 1 \\
1
\end{ytableau}
\ 
\begin{tabular}{c}
\\
\\
\\
\\
$\xrightarrow{\mathrm{flip}}$
\end{tabular}
&
\begin{ytableau}
\none&\none&\none&1\\
\none&2 & 1 &1\\
2  & {2} & {1}{2} &{1}\\
2  & 2 & 2  &1\\
2  & 2 \\
2
\end{ytableau} 
\begin{tabular}{c}
\\
\\
\\
\\
$\xrightarrow[\mathrm{res}_2]{\text{(2M)}}$
\end{tabular}
&
\begin{ytableau}
\none&\none&\none&1\\
\none&1 & 2 &1\\
{2}  & {1}2 & 2& 1\\
2  & 2 & 2  &1\\
2  & 2 \\
2
\end{ytableau} 
\begin{tabular}{c}
\\
\\
\\
\\
$\xrightarrow[\mathrm{res}_1]{\text{(2M)} }$
\end{tabular}\\

\begin{ytableau}
\none&\none&\none&1\\
\none&2 & 2 &1\\
12  & 2 & 2& 1\\
2  & 2 & 2  &1\\
2  & 2 \\
2
\end{ytableau}
\ 
 \begin{tabular}{c}
\\
\\
\\
\\
$\xrightarrow[\mathrm{res}_3]{(21)}$
\end{tabular} &

\begin{ytableau}
\none&\none&\none&2\\
\none&2 & 1 &2\\
12  & 2 & 1& 2\\
2  & 2 & 1  &2\\
2  & 2 \\
2
\end{ytableau} 
\begin{tabular}{c}
\\
\\
\\
\\
$\xrightarrow[\mathrm{res}_2]{(21)}$
\end{tabular}
&
\begin{ytableau}
\none&\none&\none&2\\
\none&1 & 2 &2\\
1 {2}  &  {1} & 2& 2\\
2  & 1 & 2  &2\\
2  & 1 \\
2
\end{ytableau} 
\begin{tabular}{c}
\\
\\
\\
\\
$\xrightarrow[\mathrm{res}_1]{\text{(M1)} }$
\end{tabular}\\
\quad
\begin{ytableau}
\none&\none&\none&2\\
\none&1 & 2 &2\\
1  & 12 & 2& 2\\
1  & 2 & 2 &2\\
1  & 2 \\
1
\end{ytableau}\begin{tabular}{c}
\\
\\
\\
\\
$=\Phi(T)$
\end{tabular} &  & 

\end{tabular}
\caption{Illustration of the map  $\Phi$.}\label{VBIO-1}
\end{figure}

\subsection{Benign 12-SV table}

For a 12-SV table $T$,  consider its mixed columns.
If   column $k$ of $T$  is mixed,
then   define $\mathrm{sep}_k
(T)$ as the smallest row index $r $ such that 
$2 \in T (r, k)$.
Suppose that 
\[
    \{k\colon \text{the $k$-th column of $T$ is mixed}\}=\{k_1<k_2<\cdots<k_p\}.
\]
Denote 
\[
\mathrm{seplist} (T)=\left(\mathrm{sep}_{k_1}(T),\mathrm{sep}_{k_2}(T),\ldots,\mathrm{sep}_{k_p}(T)\right). 
\]

\begin{definition}
We say that a 12-SV table $T$ is \textit{benign} if  
\  $\mathrm{seplist} (T)$ is a weakly decreasing sequence, with an extra  constraint  that for $j\geq 1$,    $\mathrm{sep}_{k_j}(T)>\mathrm{sep}_{k_{j+1}}(T)$ whenever there exists a (unique) box  in column $k_{j+1}$ filled with the set  $\{1,2\}$.
\end{definition}

For example,  the two 12-SV tables in Figure  \ref{fig:my_label-1}  have the same list
\[
 \mathrm{seplist} (T_1) =\mathrm{seplist} (T_2)=(3,3,2).
\]
By definition, the left one is benign, while the right one is not benign since it was supposed to satisfy  $\mathrm{sep}_{1}(T_2)>\mathrm{sep}_{3}(T_2)$ . 

\begin{figure}[ht]
    \centering
    \begin{tabular}{ccc}        
\begin{ytableau}
\none&\none&\none&1 & 2\\
\none&1  & 1  &1\bl{2}  & {2}\\
1\bl{2}   & 1 & \bl{2}  &2  & 2\\
2    & 1 & 2  &{2}  & 2\\
2    & 1 & 2\\
2
\end{ytableau}&\quad\quad &
\begin{ytableau}
\none&\none&\none&1   & 2\\
\none    &1   & 1 &\ml{2} &{2}\\
1\ml{2}  & 1  & 1\ml{2}  &2  & 2\\
2        & 1  & 2       &{2}  & 2\\
2        & 1  & 2\\
2
\end{ytableau}\\
 & & \\
$T_1$ & &$T_2$
    \end{tabular}
    \caption{$T_1$ is benign and $T_2$ is not benign.}
    \label{fig:my_label-1}
\end{figure}


Let $\mathrm{BSVT}^{12} (\lambda/\mu)$ denote the set of all benign 12-SV tables of
shape $\lambda/\mu$. Clearly, we have $\mathrm{SVRPP}^{12} (\lambda/\mu) 
\subseteq
\mathrm{BSVT}^{12} (\lambda/\mu)$.
Recall the  $\flip$ map   on 12-SV tables. For $T \in \mathrm{BSVT}^{12} (\lambda/\mu)$,    
 $\flip(T)$ still belongs to $ \mathrm{BSVT}^{12} (\lambda/\mu)$  because $T$ and $\flip(T)$ agree on   mixed
columns. So there is no obstruction to restricting the flip map to  the set    $\mathrm{BSVT}^{12} (\lambda/\mu)$.

\begin{remark}\label{RR-11}
It is easily seen that $\flip$ is an involution on $\mathrm{BSVT}^{12} (\lambda/\mu)$ that preserves all the statistics $\ceq$, $\ex$ 
 and $\mathrm{seplist}$, but switches the first two entries of the statistic $\ircont$.    
\end{remark}

\subsection{Descent resolution}

 Let $T \in \mathrm{BSVT}^{12} (\lambda/\mu)$. Assume that
$k$ is a descent of $T$. Then  the $k$-th column of $T$ is either 2-pure or mixed, and the $(k+1)$-th
column of $T$ is either 1-pure or mixed. Note that  the $k$-th and the $(k+1)$-th columns
of $T$ cannot be both  mixed because   $T$ is benign.  Therefore, columns $k$ and $k+1$ of $T$  are possibly of type (M1), or (2M), or (21), as  aforementioned in  Subsection \ref{sec:method}. After applying the descent-resolution operation $\mathrm{res}_k$ in   Figure \ref{fig:des-res-preliminary}  to $T$, the resulting filling  $\mathrm{res}_k(T)$ is still a benign 12-SV table since
\[\mathrm{seplist}(\mathrm{res}_k(T)) = \mathrm{seplist} (T).\]




The following statements can be  verified  directly from the  definitions of the $\mathrm{flip}$ map and the descent-resolution operation $\mathrm{res}_k$.

\begin{proposition}\label{HGP-25}
Assume that  $T \in \mathrm{BSVT}^{12} (\lambda/\mu)$, and $k  $ is a descent
of $T$. 
\begin{itemize}
    \item[(1)] 
    The operation $\mathrm{res}_k$ preserves all the statistics $\ceq$, $\ex$,
 and $\mathrm{ircont}$, that is, 
$$\ceq (\mathrm{res}_k(T)) = \ceq (T),\ \ex(\mathrm{res}_k(T)) =\ex(T),  \ \ircont(\mathrm{res}_k(T)) = \ircont(T).$$

\item[(2)] The integer $k$ is a descent of $\mathrm{flip} (\mathrm{res}_k (T))$, and moreover  we have
\[
 \mathrm{res}_k (\mathrm{flip} (\mathrm{res}_k (T))) = \mathrm{flip} (T).   
\]

\item[(3)]  We have
$\ell(T) > \ell(\mathrm{res}_k (T))$. Here
    $$\ell(T)=\sum_{j\geq 1}j\cdot \mathrm{sig}(\text{the j-th column of T}),$$
where, for a column $C$ of $T$, 
\begin{align*}
         \mathrm{sig}(C)=\begin{cases}
 0, & \text{ if $C$ is empty or 2-pure,}  \\
 1, & \text{ if $C$ is mixed, } 
 \\
 2, & \text{ if $C$ is 1-pure}.
\end{cases}
     \end{align*}
\end{itemize}
\end{proposition}

\begin{myproof}
The assertions in (1), (2) and (3)   in the non-set-valued case have been checked  in  \cite{Galashin 1}, that is, in Figure  \ref{fig:des-res-preliminary}, both  $*$'s in (M1) are    equal to $\{2\}$, and  both  $\star$'s in (2M) are   equal to $\{1\}$. It is not hard to verify these assertions when   every   $*$  in (M1) or every   $\star$  in (2M) is   equal to the set $\{1,2\}$.
\end{myproof}

\subsection{The descent-resolution relation}
 
Based on the descent resolution, we  define a binary relation $\rightarrow$ on  benign 12-SV tables.
Let $P, Q \in \mathrm{BSVT}^{12} (\lambda/\mu)$. Denote
$P \xrightarrow{k} Q$ if  $Q = \mathrm{res}_k (P)$ where $k$  is a descent $P$. Write $P \rightarrow Q$ if   $P \xrightarrow{k} Q$ for some descent $k$ of $P$.

 The following is just a restatement of Proposition \ref{HGP-25}.

\begin{proposition}\label{HGP-27}
Assume that  $P,Q \in \mathrm{BSVT}^{12} (\lambda/\mu)$ and  $P \rightarrow Q$. Then we have 
\begin{itemize}
    \item[(1)]   $\ceq (Q) = \ceq (P)$,  $\ex(Q)=\ex(P)$, and $\ircont(Q) = \ircont(P) $;

  \item[(2)]  
$\mathrm{flip} (Q) \rightarrow \mathrm{flip} (P)$;

    \item[(3)]   $\ell(P) >\ell(Q)$.
\end{itemize}
\end{proposition}

We now define $\xrightarrow{*}$ to be the  transitive closure of the above relation. That is,  for $P, Q \in \mathrm{BSVT}^{12} (\lambda/\mu)$, write  $P \xrightarrow{*} Q$ if   there exists a sequence \[P=P_0\rightarrow P_1\rightarrow \cdots \rightarrow P_m=Q,\ \ \ \text{where $m\geq 0$}.\]  
By Proposition \ref{HGP-27}, we immediately have 

\begin{proposition}\label{HGP-28}
Assume that  $P,Q \in \mathrm{BSVT}^{12} (\lambda/\mu)$ and  $P \xrightarrow{*} Q$. Then we have 
\begin{itemize}
    \item[(1)]   $\ceq (Q) = \ceq (P)$,  $\ex(Q)=\ex(P)$, and $\ircont(Q) = \ircont(P) $;

  \item[(2)]  
$\mathrm{flip} (Q) \xrightarrow{*} \mathrm{flip} (P)$;

    \item[(3)]   $\ell(P) \geq \ell(Q)$ (and $\ell(P) > \ell(Q)$ if $P\neq Q$).
\end{itemize}
\end{proposition}

The following lemma will be used    in the proof of Theorem  \ref{HGT-22}.

\begin{lemma}\label{HGL-29}
Assume that  $T, P, Q \in \mathrm{BSVT}^{12} (\lambda/\mu)$ such that $T \rightarrow P$
and $T \rightarrow Q$. Then, there exists  $T' \in \mathrm{BSVT}^{12} (\lambda/\mu)$ such that $P\xrightarrow{*}T'
$ and $Q\xrightarrow{*}T'$.    
\end{lemma}

\begin{myproof}
When $P = Q$, we may simply  let $T'=P=Q$. Next, assume that $P\neq Q$. Suppose that $T\xrightarrow{i}P$ and $T \xrightarrow{j}
Q$, namely, $P
 = \mathrm{res}_i (T)$ and
$Q = \mathrm{res}_j (T)$. Without loss of generality,   let  
$i < j$. The arguments are divided into two cases.

Case 1: $i < j - 1$. 
It is clear  that $\mathrm{res}_j (\mathrm{res}_i (T))=\mathrm{res}_i (\mathrm{res}_j (T))$. Then  $T' = \mathrm{res}_i (\mathrm{res}_j (T))$ is a desired benign 12-SV table.

Case 2:  $i = j - 1$. 
In this case,   both the numbers $1$ and $2$ appear in the $j$-th column of $T$. So this column is mixed. Together with the fact that  $T$ is a benign 12-SV table such that both  $j - 1$ and $j$ are descents, we obtain that the $(j - 1)$-th column of $T$ is 2-pure and the $(j + 1)$-th column of $P$
is 1-pure. In view of the descent resolutions in  Figure \ref{fig:des-res-preliminary}, 
it is easy to check that  
$$\mathrm{res}_{j-1} (\mathrm{res}_j(\mathrm{res}_{j-1}(T))) = \mathrm{res}_j (\mathrm{res}_{j-1}(\mathrm{res}_j (T))),$$
which means $ \mathrm{res}_{j-1} (\mathrm{res}_j (P)) = \mathrm{res}_j (\mathrm{res}_{j-1} (Q))$.  Let  $T'$ be  
this 12-SV table.  
\end{myproof}

\subsection{Proof of Theorem  \ref{HGT-22}}

\begin{proposition}\label{HGP-210}
 For each $T \in \mathrm{BSVT}^{12} (\lambda/\mu)$, there exists a unique $N$ belonging to $  
\sr^{12} (\lambda/\mu)$ such  that $T\xrightarrow{*}N$  
\end{proposition}

\begin{myproof}
Suppose that  $T \in \mathrm{BSVT}^{12} (\lambda/\mu)$.  Set
$$\mathrm{Norm} (T)=\{N \in \sr^{12} (\lambda/\mu)\colon T\xrightarrow{*}N\}.$$
We   show that  $\mathrm{Norm} (T)$ contains exactly one element.
The argument is by  induction on $\ell(T)$. The initial condition is when $\ell(T)=0$. In this case, it is trivial that $\mathrm{Norm} (T)=\{T\}$ since $T$ is a SVRPP with each column being 2-pure. We next consider the case of $\ell(T)>0$. 
Suppose  that  {$\# \mathrm{Norm} (S)=1$  for every $S \in \mathrm{BSVT}^{12} (\lambda/\mu)$ with $\ell (S) < \ell (T)$.}

Let 
\[
  \mathbf{Z} =\{Z \in \mathrm{BSVT}^{12} (\lambda/\mu) \colon T \rightarrow Z
\}.  
\]
be the set of all benign
12-SV tables  which can be obtained from $T$ by resolving one descent.  
If $\textbf{Z}$ is empty,
then $T \in
\sr^{12} (\lambda/\mu)$, yielding that $\mathrm{Norm} (T) = \{T\}$. It remains to consider the situation  when $\textbf{Z}$ is nonempty. Note that $T \notin
\sr^{12} (\lambda/\mu)$ when $\textbf{Z}$ is nonempty.

Keep in mind that for every $N \in
\sr^{12} (\lambda/\mu)$, if  $ T
 \xrightarrow{*}N$, then $Z
\xrightarrow{*} N$ for some
$Z \in \textbf{Z}$. This implies that every $ N \in \mathrm{Norm} (T)$ must belong to $\mathrm{Norm} (Z)$ for some
$Z \in \textbf{Z}$. Therefore,
\begin{equation}\label{BgYt-3}
\mathrm{Norm} (T) =\bigcup_{Z \in \textbf{Z}}
\mathrm{Norm} (Z).
\end{equation}
We have the following observations.
\begin{itemize}
    \item[(I)]  By (3) in  Proposition \ref{HGP-27}, we see that $\ell(T)>\ell(Z)$ for any $Z \in \textbf{Z}$. By induction,   $\mathrm{Norm} (Z)$ is a one-element set.

    \item[(II)] For   $P, Q \in \textbf{Z}$, we have $\mathrm{Norm} (P)\cap \mathrm{Norm} (Q) \neq \emptyset$. This can be seen as follows. According to Lemma \ref{HGL-29}, there exists $T'\in \mathrm{BSVT}^{12} (\lambda/\mu)$ such that $P\xrightarrow{*}T'
$ and $Q\xrightarrow{*}T'$. By Propositions \ref{HGP-27} and \ref{HGP-28}, we see that $\ell(T) > \ell (P) \geq \ell (T')$. By induction,   $\mathrm{Norm} (T')$ has one element. It is evident that $\mathrm{Norm} (T')\subseteq \mathrm{Norm} (P)\cap \mathrm{Norm} (Q)$, and so $\mathrm{Norm} (P)\cap \mathrm{Norm} (Q) \neq \emptyset$. 
\end{itemize}

By (I) and (II),   all the one-element sets in the union \eqref{BgYt-3} are identical.  Hence  $\mathrm{Norm} (T)$  itself has exactly one element.  
\end{myproof}

For $T \in \mathrm{BSVT}^{12} (\lambda/\mu)$, by Proposition \ref{HGP-210},  we define $\mathrm{norm} (T)$ to be the unique  $N \in
\sr^{2} (\lambda/\mu)$ such that $T\xrightarrow{*}N$.

We are now in a position to give a proof of Theorem  \ref{HGT-22}.

\noindent 
{\it Proof of Theorem  \ref{HGT-22}.}
We first show that the map ${\Phi}$ defined in Subsection \ref{sec:method} is a well-defined involution. Let $T\in \sr^{12}(\lambda/\mu)$.  In the process of  applying the descent-resolution operations    to $\mathrm{flip}(T)$, according to (3) in Proposition  \ref{HGP-27} and Proposition \ref{HGP-210}, it eventually terminates at  ${\Phi} (T) = \mathrm{norm} (\mathrm{flip} (T))$. Since \[\mathrm{flip} (T) \xrightarrow{*} \mathrm{norm} (\mathrm{flip} (T)) =
{\Phi} (T),\] 
it follows from (2) in Proposition \ref{HGP-28} that 
\begin{equation}\label{fjhI-2}
  \mathrm{flip} ({\Phi} (T)) \xrightarrow{*} \mathrm{flip} (\mathrm{flip} (T)) = T.  
\end{equation}
On the other hand,  ${\Phi}({\Phi}(T)) = \mathrm{norm} (\mathrm{flip} ({\Phi} (T)))$ is the unique $ N \in \sr^{12}(\lambda/\mu)$ such that $\mathrm{flip} ({\Phi} (T)) \xrightarrow{*} N$. 
This, together with \eqref{fjhI-2},
forces that   ${\Phi}({\Phi}(T)) = T$, concluding  that $\Phi$ is  an involution. 
 
Because $\mathrm{flip} (T) \xrightarrow{*}  {\Phi} (T)$, by Remark \ref{RR-11} and (1) in Proposition   \ref{HGP-28},  $\Phi$
preserves   the statistics $\ceq$ and $\ex$, but switches the first two entries of the statistic $\ircont$. So the proof is complete. 
\qed

\section{Proof of Theorem \ref{THM-2}}\label{Sc-43}

In this section, we build a $\gl_n$-crystal structure on the set $\mathrm{SVRPP}^n(\lm)$.  We prove that the character of each connected component in its crystal graph is equal to a Schur polynomial.  As an application, we obtain the expansion of ${G}_{\lm}(\mathbf{x}_n ;\mathbf{t};\mathbf{w})$ in  terms of  Schur polynomials, finishing   the  proof of 
Theorem \ref{THM-2}. 

\subsection{\texorpdfstring{$\gl_n$-crystal}{gl(n)-crystal}}

We mostly follow the terminology and notation in \cite{Bump,MTY}.
Let $n$ be a given positive integer. Assume that  $\mathcal{B}$ is a set with maps $\mathrm{wt}\colon   \mathcal{B}\to \ZZ^n$
and 
$e_i,f_i\colon   \mathcal{B} \to \mathcal{B} \sqcup \{0\}$  for $i \in [n-1]$,
where $0 \notin \mathcal{B}$.
The map $\mathrm{wt}$ is called  the  {\textit{weight map}}, and   $e_i$ and $f_i$ are called   {\textit{raising}} and {\it {lowering operators}}, respectively.

Write $\e_1,\e_2,\dots,\e_n\in\ZZ^n$ for the standard basis.
The set $\mathcal{B}$ is a (Kashiwara)
  {$\gl_n$-{\it crystal}}
if
for every $i \in [n-1]$ and $b,c \in \mathcal{B}$,
   $e_i(b) = c $ if and only if $ f_i(c) = b$,
in which case 
\[\mathrm{wt}(c) = \mathrm{wt}(b) + \e_{i} -\e_{i+1}.\]
The  {\textit{character}} is the weight generating function in $\mathbf{x}_n=(x_1,\ldots, x_n)$:
\[ \mathrm{ch}(\mathcal{B}) := \sum_{b \in \mathcal{B}} \mathbf{x}_n^{\mathrm{wt}(b)}.\]

Define the  {string lengths} $\varepsilon_i, \varphi_i\colon  \mathcal{B} \to \{0,1,2,\dots\}\sqcup \{\infty\}$
by 
$$\varepsilon_i(b) := \max\left\{ k\geq 0 \colon e_i^k(b) \neq 0\right\},\ \ \ \ 
\varphi_i(b) := \max\left\{ k \geq 0\colon  f_i^k(b) \neq 0\right\}.
$$
A $\gl_n$-crystal $\mathcal{B}$ is  {\textit{seminormal}}
if the string lengths   take finite values
such that  for any $b \in \mathcal{B}$,
\[\varphi_i(b) - \varepsilon_i(b) = \mathrm{wt}(b)_i - \mathrm{wt}(b)_{i+1}.
\]
where $\mathrm{wt}(b)_i$ is the $i-$th component of $\mathrm{wt}(b)$.

The associated {\it crystal graph} of  $\mathcal{B}$ is a directed graph  with vertex set $\mathcal{B}$, such that  $b \xrightarrow{i} c$ is a (labeled) edge if and only if $ f_i(b)=c$. 
Each connected component of the crystal graph of $\mathcal{B}$ inherits  a crystal structure from  $\mathcal{B}$, which is  called  a {\textit{full subcrystal}}. An element $b\in \mathcal{B}$  is called a {\textit{highest weight element}}  if it is a source vertex in the   crystal graph (that is, $e_i(b)=0$ for all $i\in [n-1]$). 
Let $\HW(\mathcal{\mathcal{B}})$ denote the set of highest weight elements in $\mathcal{B}$.
Two crystals  are {\it isomorphic}  if there is a weight-preserving graph isomorphism between their
crystal graphs.

A $\gl_n$-crystal is  called {\textit{normal}} if it is isomorphic to a disjoint union of full subcrystals in the tensor powers of the {\it standard} $\gl_n$-crystal $\mathbb{B}_n$. Here the crystal graph of $\mathbb{B}_n$ is   the path on $[n]$ with edges $i \xrightarrow{i} i+1$ and weight map $\mathrm{wt}(i)=\mathbf{e}_i$. Each full subcrystal in the tensor powers of $\mathbb{B}_n$ is indexed by a partition $\lambda=(\lambda_1\geq \lambda_2\geq \cdots \geq \lambda_n\geq 0)$, and   will be denoted  by $\mathcal{B}(\lambda)$. 



The crystal  $\mathcal{B}(\lambda)$ corresponds to the crystal basis of the highest module $V (\lambda)$ of the quantum group $U_q(\gl_n)$ \cite{Kas1,Kas2}.   
This crystal  can be constructed on semistandard tableaux of shape $\lambda$,  with crystal operators defined relying  on the paring rule as defined in Subsection \ref{IHN_09}. 
Here we shall use an isomorphic realization based on words in the alphabet $[n]$. We refer to  \cite{Galashin 2} for a proof of the equivalence between these two constructions. 

Let $S_{n}^m:=[n]^m$ be the set of all words of length $m$ with letters in $[n]$. For a word $s=(s_1,\dots,s_m)\in S_{n}^m$, the weight $\mathrm{wt}(s)$ is defined as the vector $(c_1,\dots,c_n)$ where $c_i$ is the number of occurrences of $i$ in $s$. 
For $i\in [n-1]$,   define the operators $E_i,  F_i:S_{n}^m\to S_{n}^m\cup\{0\}$ using the paring rule.
To be specific, for $s=(s_1,\dots, s_m)\in S_{n}^m$,
consider the subword  of $s$ occupied by   $i$ or $i + 1$.
Replace each $i$ by $\bplus$, and each $i+1$ by $\bminus$. Next, cancel successively  the matched  ordered $(\bminus, \bplus)$-pairs. This results in  a string  of the form $\bplus \cdots \bplus \bminus \cdots \bminus$.
Define 
\begin{itemize}
    \item $E_i(s)$: If there is no $\bminus$ in the resulting string, then set $E_i(s)= 0$. Otherwise, $E_i(s)$ is obtained from $s$ by replacing  the $i+1$ corresponding to  the leftmost unpaired  $\bminus $ by an $i$.
    
    \item $F_i(s)$: If there is no $\bplus$ in the resulting string, then set $F_i(s)= 0$. Otherwise, $F_i$ is obtained from $s$ by replacing  the $i$ corresponding to  the rightmost unpaired  $\bplus $ by an $i+1$.
\end{itemize}

Together with the above operators,  
the set
$S_{n}^m$  forms a normal $\gl_n$-crystal structure. 
 We collect some of its properties. 

\begin{proposition}[\cite{Galashin 2}]
\label{lemma:crystal}
\begin{itemize}
 \item[(1)]  Each connected component (namely, full subcrystal) in the   crystal graph of $S_{n}^m$ has exactly one highest weight element;
 
 \item[(2)] Each connected component is completely determined (up to  isomorphism) by the weight  of the unique highest weight element, which is   a partition $\lambda$;
 
 \item[(3)] The character of the connected component corresponding to $\lambda$ is equal to the Schur polynomial $s_{\lambda}(\mathbf{x}_n)$.  
\end{itemize}
\end{proposition}
 
Now, the crystal $\mathcal{B}(\lambda)$ is isomorphic to the full subcrystal  of $S_{n}^m$ determined by the partition $\lambda$.

\subsection{Crystal structure on set-valued reverse plane partition}\label{IHN_09}

We begin  by  defining   the crystal operators  on $\mathrm{SVRPP}^n(\lm)$.
Let $T\in \mathrm{SVRPP}^n(\lm)$. Define the weight map as \[\mathrm{wt}(T)=\mathrm{ircont}(T).\]

Fix $i\in [n-1]$. To define the  operators $e_i$ and   $f_i$, recall that  $T^{(i)}$ is the subtableau of  $T$ occupied by $i$ or $i+1$. Imitating the description in  Section \ref{sec:method}, a column of $T^{(i)}$ is called $i$-pure (resp., $(i+1)$-pure) if it contains only $\{i\}$'s (resp., $\{i+1\}$'s), and otherwise it is called mixed. 

\begin{figure}[ht]
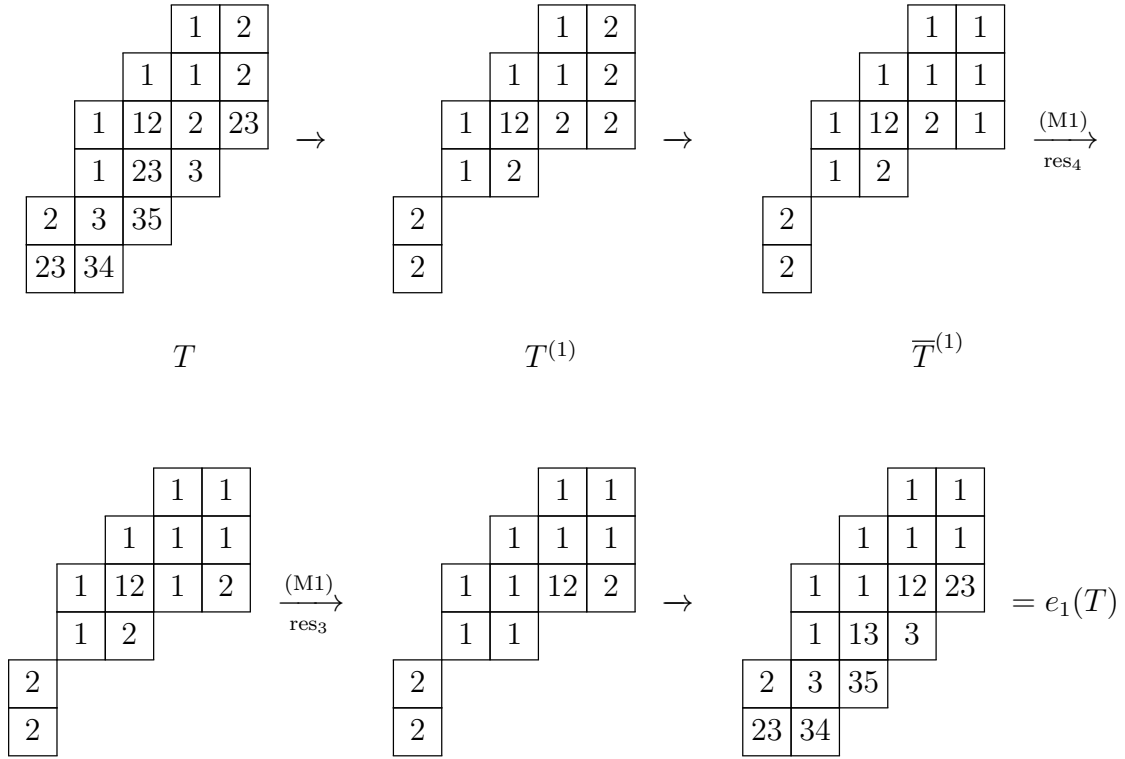


\begin{tabular}{ccc}
\begin{ytableau}
\none & \none& \none &  1    & 2   \\
\none & \none &  1    &  1    & 2    \\
\none &  1    &  12   &  2    & 23    \\
\none &  1   &  23   &  3    \\
 2    &  3    &  35 \\
 23    &  34    
\end{ytableau}
\begin{tabular}{c}
\\
\\
\\
\\
\\
$\xrightarrow{}$
\end{tabular}
&
\begin{ytableau}
\none & \none& \none &  1    & {2}   \\
\none & \none &  1    &  1    & {2}    \\
\none &  1     &  12   &  2    & {2}    \\
\none &  1    &  2    \\
 2     \\
 2         
\end{ytableau}
\begin{tabular}{c}
\\
\\
\\
\\
\\
$\xrightarrow{}$
\end{tabular}
&
\begin{ytableau}
\none & \none& \none &  1    & 1   \\
\none & \none&  1    &  1    & 1    \\
\none &  1    &  12   &  {2}    & {1}    \\
\none &  1    &  2    \\
 2     \\
 2         
\end{ytableau}
\begin{tabular}{c}
\\
\\
\\
\\
\\
$\xrightarrow[\mathrm{res}_4]{\mathrm{(M1)}}$
\end{tabular}\\
& & \\
$ T$& $T^{(1)}$ & $\overline{T}^{(1)}$\\
\begin{ytableau}
\none & \none& \none &  1    & 1   \\
\none & \none&  1    &  1    & 1    \\
\none &  1   &  1{2}  & {1}  & 2  \\
\none &  1   &  2    \\
 2     \\
 2         
\end{ytableau}
\begin{tabular}{c}
\\
\\
\\
\\
\\
$\xrightarrow[\mathrm{res}_3]{\mathrm{(M1)}}$
\end{tabular}&
\begin{ytableau}
\none & \none& \none &  1    & 1   \\
\none & \none&  1    &  1    & 1    \\
\none &  1   &  1   &  12    & 2    \\
\none &  1   &  1    \\
 2     \\
 2         
\end{ytableau}
\begin{tabular}{c}
\\
\\
\\
\\
\\
$\xrightarrow{}$
\end{tabular}&
\begin{ytableau}
\none & \none& \none &  1    & 1   \\
\none & \none &  1    &  1    & 1    \\
\none &  1     &  1   &  12    & 23    \\
\none &  1     &  13  &  3  \\
 2    &  3     &  35 \\ 
 23   &  34       
\end{ytableau}
\begin{tabular}{c}
\\
\\
\\
\\
\\
$=e_1(T)$
\end{tabular}\\
\end{tabular}

\vspace{0.5cm}
\caption{Illustration of the crystal operator $e_i$ with $i=1$.}\label{U9I-1}
\end{figure}

\begin{figure}[ht]
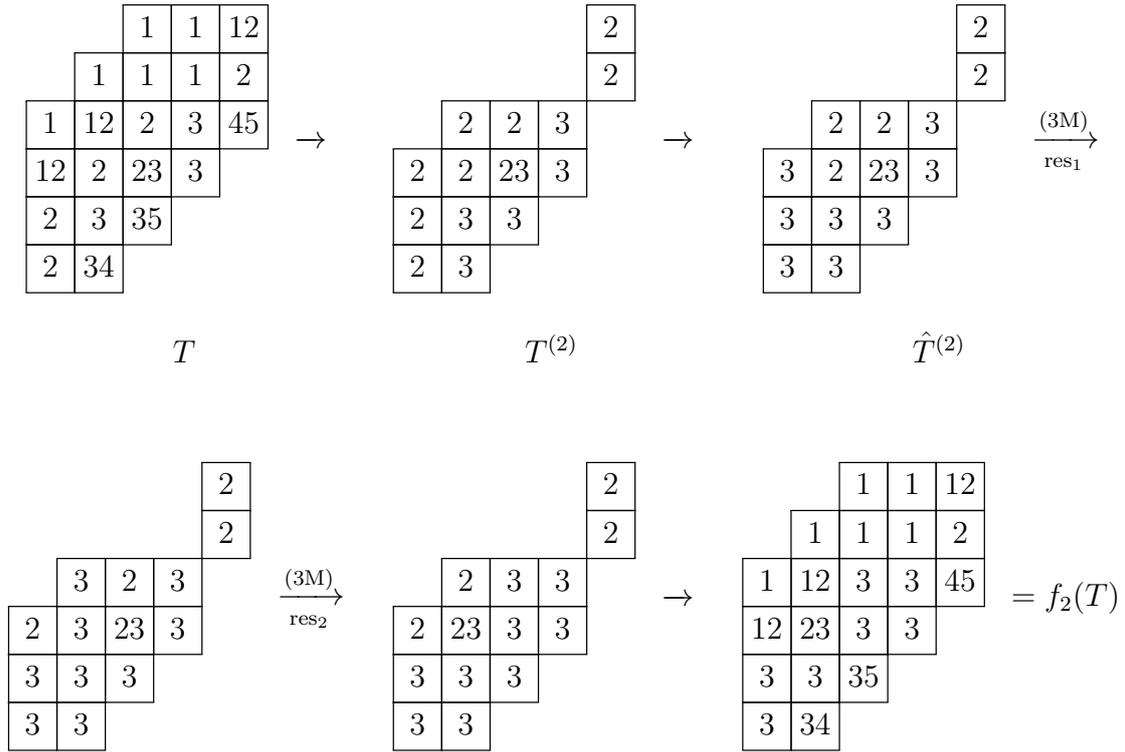


\begin{tabular}{ccc}
\begin{ytableau}
\none &\none& 1 &  1    & 12   \\
\none & 1 &  1    &  1    & 2    \\
 1 &  12    &  2   & 3    &45    \\
 12 &  2   &  23   &  3    \\
 2    &  3    &  35 \\
 2    &  34    
\end{ytableau}
\begin{tabular}{c}
\\
\\
\\
\\
\\
$\xrightarrow{}$
\end{tabular}&
\begin{ytableau}
\none &\none& \none &  \none    & 2   \\
\none & \none & \none   &  \none   & 2    \\
 \none &  2    &  2   & 3        \\
 {2} &  2   &  23   &  3    \\
 {2}  &  3    &  3 \\
{2}  &  3    
\end{ytableau}
\begin{tabular}{c}
\\
\\
\\
\\
\\
$\xrightarrow{}$
\end{tabular}&
\begin{ytableau}
\none &\none& \none &  \none    & 2   \\
\none & \none & \none   &  \none   & 2    \\
 \none &  2    &  2   & 3        \\
{3} &  {2}   &  23   &  3    \\
 3    &  3    &  3 \\
 3    &  3    
\end{ytableau}
\begin{tabular}{c}
\\
\\
\\
\\
\\
$\xrightarrow[\mathrm{res}_1]{\mathrm{(3M)}}$
\end{tabular}\\
& & \\
$T$& $T^{(2)}$ &$\hat{T}^{(2)}$ \\
\begin{ytableau}
\none &\none& \none &  \none    & 2   \\
\none & \none & \none   &  \none   & 2    \\
 \none &  3    &  2   & 3        \\
 2     &  {3}   &  {2}3   &  3    \\
 3    &  3    &  3 \\
 3    &  3    
\end{ytableau}
\begin{tabular}{c}
\\
\\
\\
\\
\\
$\xrightarrow[\mathrm{res}_2]{\mathrm{(3M)}}$
\end{tabular}
&
\begin{ytableau}
\none &\none& \none &  \none    & 2   \\
\none & \none & \none   &  \none   & 2    \\
 \none &  2    &  3   & 3        \\
 2     &  23   &  3   &  3    \\
 3    &  3    &  3 \\
 3    &  3    
\end{ytableau}
\begin{tabular}{c}
\\
\\
\\
\\
\\
$\xrightarrow{}$
\end{tabular}
&
\begin{ytableau}
\none &\none& 1 &  1    & 12   \\
\none & 1 &  1    &  1    & 2    \\
 1 &  12    &  3   & 3    &45    \\
 12 &  2 3  &  3   &  3    \\
 3    &  3    &  35 \\
 3    &  34    
\end{ytableau}
\begin{tabular}{c}
\\
\\
\\
\\
\\
$=f_2(T)$
\end{tabular}\\
& & \\
\end{tabular}

\caption{Illustration of crystal operator $f_i$ with $i=2$.}\label{U9I-2}
\end{figure}

Put a `$\bplus$' above each $i$-pure column of $T^{(i)}$, and a  `$\bminus$' above each $(i+1)$-pure column of $T^{(i)}$ (ignore the mixed columns). 
This results in a string  on the alphabet $\{\bplus, \bminus\}$, which will be called the $i$-{\it signature} of $T$.
Next, cancel successively  the matched  ordered $(\bminus, \bplus)$-pairs  until yielding a string  of the form $\bplus \cdots \bplus \bminus \cdots \bminus$, which is  called the \textit{reduced} $i$-signature of $T$.
\begin{itemize}
\item[$\bullet$] $e_i(T)$: If there is no $\bminus$ in the resulting string, then set $e_i (T) = 0$. Otherwise, let $A$ be  the $(i+1)$-pure  column of $T^{(i)}$ corresponding to the leftmost unpaired $\bminus$. Replace all the $\{i+1\}$'s in $A$ by $\{i\}$'s, and the resulting filling is denoted $\overline{T}^{(i)}$. 
Clearly, $\overline{T}^{(i)}$ may have a descent occurring at the position of the column  immediately to the left of $A$. 
Note that $\overline{T}^{(i)}$ is   benign. 
Let $e_i(T)$ be the unique SVRPP obtained from $T$ by performing  the descent-resolution algorithm  to $\overline{T}^{(i)}$ (we treat $i$ and $i+1$ as $1$ and $2$, respectively). In other words, $e_i(T)$ is obtained from $T$ by replacing ${T}^{(i)}$ with $\mathrm{norm}(\overline{T}^{(i)})$. See Figure  \ref{U9I-1}  for an illustration of  the construction   of $e_i(T)$ for $i=1$.

\item[$\bullet$] $f_i(T)$: If there is no $\bplus$ in the resulting string, then set $f_i (T) = 0$. Otherwise, let $A$ be the $i$-pure column of of $T^{(i)}$ corresponding to the rightmost unpaired  $\bplus$. Replace all the $\{i\}$'s in $A$ by $\{i+1\}$'s, and the result filling is denoted $\hat{T}^{(i)}$. Let $f_i(T)$ be the unique SVRPP obtained from $T$ by performing  the descent-resolution algorithm to  $\hat{T}^{(i)}$.  See Figure  \ref{U9I-2}  for an illustration of  the construction   of $f_i(T)$ for $i=2$.
\end{itemize}

\begin{remark}
The crystal  operators on $\mathrm{SVRPP}^n(\lm)$ generalize simultaneously   the   crystal operators on set-valued tableaux  by Monical, Pechenik and Scrimshaw \cite{Monical} and on  reverse plane partitions by Galashin  \cite{Galashin 2}.   
\end{remark}

We now show  that    the above operators  satisfy the crystal axioms.  

\begin{proposition}\label{kou98}
    Let $T, T' \in \mathrm{SVRPP}^n(\lm)$. Then
    \begin{itemize}
        \item[(1)]$e_i (T), f_i (T) \in \mathrm{SVRPP}^n(\lm) \sqcup \{0\}$;
        
        \item[(2)]$e_i (T) = T' \Longleftrightarrow T = f_i (T')$;
        
\item[(3)]
        If 
   $e_i(T) = T' $, then  
\[\mathrm{wt}(T') = \mathrm{wt}(T) + \e_{i} -\e_{i+1}.\]
    \end{itemize}

\end{proposition}

\begin{myproof}
(1) This follows from the definitions of $e_i$ and $f_i$.

(2) Let first check that $e_i (T) = T'$ implies $ T = f_i (T')$. As before,  use $T^{(i)}$ to denote  the subtableau of  $T$ occupied by $i$ or $i+1$. We shall analyze  how the desent-resolution algorithm is implemented in the construction of $e_i(T)$.  

Let $A_i$ be the $(i+1)$-pure column of $T^{(i)}$ marked with the leftmost unpaired  $\bminus$,  $B_i$ be the leftmost $(i+1)$-pure column of $T^{(i)}$ lying 
to the right of $A_i$, and  $C_i$ be the rightmost $i$-pure column lying to the left of $A_i$. By the above choices, it is easy to check that there
is no $i$-pure column   between $A_i$ and $B_i$, and no $(i+1)$-pure column  between
$C_i$ and $A_i$. 
In other words,  all the columns of $T^{(i)}$ lying between $C_i$ and
$A_i$, as well as between $A_i$ and $B_i$, are mixed. Of course, it is possible that  either $B_i$ or $C_i$ does not exist, and in such case,
all the columns to the right (resp., to the left) of $A_i$ are mixed.

 We explain  that only the descent-resolution of type  (M1)$\rightarrow$(1M) in Figure  \ref{fig:des-res-preliminary} will appear when applying  $e_i$ to get $T'$. Since the column $A_i$ becomes an $i$-pure column,  the only possible descent could occur at the column $A_i-1$ (which is the column immediately  to the left of $A_i$),
and  as a result, the $i$-pure column switches  to the left. Now this $i$-pure column is surrounded by mixed columns. So the only possible descent resolution  
could occur  between this $i$-pure column and the mixed column to its left. Continue the (M1) type descent resolution (if any)  until there exists no descent. Note that  the $i$-pure column in the descent-resolution procedure can  move leftwards at most to  the position of the column $C_i+1$ since  $C_i$ is   $i$-pure.

Let $D_i$ be the column where the  descent-resolution procedure  stops. Clearly, we have $C_i<D_i\leq A_i$. By the above analysis,   $D_i$ is exactly  the column corresponding to  the rightmost uncanceled $\bplus$ in $(T')^{(i)}$. To construct $f_i(T')$, we need only to reverse   the above descent-resolution steps (which are of type (2M)$\rightarrow$(M2)) between $D_i$ and $A_i$. As a result, this will recover $T$, implying that  $f_i (T')=T$.
Using completely similar arguments, we can check   that  if $T = f_i (T')$, then  $e_i (T) = T'$.

(3) This is clear from the definition of $e_i$. 
\end{myproof}


Recall that $\varepsilon_i(T)$ (resp., $\varphi_i(T)$) is the largest $k$ such that $e_i^k(T) \neq 0$ (resp., $f_i^k(T) \neq 0$). 
The proof of Proposition \ref{kou98} leads to the following observation. 

\begin{corollary}\label{cnvn-3}
For $T\in \mathrm{SVRPP}^n(\lm)$,   $\varepsilon_i(T)$ (resp.,  $\varphi_i(T)$) equals the number of   minus (resp., plus) signs in the reduced $i$-signature of $T$.
\end{corollary}

Proposition \ref{kou98} and Corollary  \ref{cnvn-3} together imply that $\mathrm{SVRPP}^n(\lm)$ is a seminomal $\gl_n$-crystal.

\begin{corollary}\label{gsgl-3}
With the crystal operators $\mathrm{wt}$, $e_i$ and $f_i$ defined above, $\mathrm{SVRPP}^n(\lm)$ is a  seminormal $\gl_n$-crystal. 
\end{corollary}
 
\begin{myproof}
Proposition \ref{kou98} says that $\mathrm{SVRPP}^n(\lm)$ is a $\gl_n$-crystal. 
By Corollary  \ref{cnvn-3}, it follows that for $T\in \mathrm{SVRPP}^n(\lm)$, $\varphi_i(T) - \varepsilon_i(T)$  is equal to  the value $ \mathrm{wt}(T)_i - \mathrm{wt}(T)_{i+1}$. So this crystal is seminormal. 
\end{myproof}

\begin{remark}
If a  crystal $\mathcal{B}$  is seminormal, then its character $\mathrm{ch}(\mathcal{B})$ is symmetric   \cite{Bump}. So Corollary  \ref{gsgl-3} also suggests the symmetry of $G_{\lambda/\mu}(\mathbf{x}_n;\mathbf{1}; \mathbf{1})$. 
\end{remark}

We may now consider the crystal graph of  $\mathrm{SVRPP}^n(\lm)$.
In Figure \ref{dnk-2-1}, we list two of the connected components in the crystal graph of $\sr^n(\lm)$ with $\l=(3,2)$, $\m=(1)$ and $n=3$. Moreover, 
Appendix \ref{hgi-6-9} gives a full list of  connected components of the crystal graph for $\lambda=(2,2)$, $\mu=(1)$ and $n=3$.
\begin{figure}[h t]
\[
\ytableausetup{boxsize=2em}
\scalebox{0.7}{
\begin{tikzpicture}[>=latex]

\node (A) {\ytableaushort{\none {1}{1},{1}{12}}};
\node[below  = 1  of A] (B) 
{\ytableaushort{\none {1}{2},{1}{12}}};
\node[below right = 1 and 0.5 of A] (C) 
{\ytableaushort{\none {1}{1},{1}{13}}};

\node[below left = 1 and 0.5 of B] (D) 
{\ytableaushort{\none {2}{2},{12}{2}}};
\node[below = 1  of B] (E) 
{\ytableaushort{\none {1}{3},{1}{12}}};
\node[below = 1  of C] (F) 
{\ytableaushort{\none {1}{2},{1}{13}}};

\node[below = 1  of D] (G) 
{\ytableaushort{\none {2}{3},{12}{2}}};
\node[below = 1  of E] (H) 
{\ytableaushort{\none {1}{3},{1}{13}}};
\node[below = 1  of F] (I) 
{\ytableaushort{\none {2}{2},{1}{23}}};

\node[below = 1  of G] (J) 
{\ytableaushort{\none {3}{3},{12}{3}}};
\node[below = 1  of H] (K) 
{\ytableaushort{\none {2}{3},{1}{23}}};
\node[below = 1  of I] (L) 
{\ytableaushort{\none {2}{2},{2}{23}}};

\node[below = 1  of J] (M) 
{\ytableaushort{\none {3}{3},{13}{3}}};
\node[below = 1  of K] (N) 
{\ytableaushort{\none {2}{3},{2}{23}}};
\node[below = 1  of N] (O) 
{\ytableaushort{\none {3}{3},{2}{23}}};

\path (A) edge[pil,color=blue] node[ left,black]{$1$} (B);
\path (A) edge[pil,color=red] node[above right,black]{$2$} (C);

\path (B) edge[pil,color=blue] node[above  left,black]{$1$} (D);
\path (B) edge[pil,color=red] node[right,black]{$2$} (E);
\path (C) edge[pil,color=blue] node[ right,black]{$1$} (F);

\path (D) edge[pil,color=red] node[left,black]{$2$} (G);
\path (E) edge[pil,color=blue] node[above left,black]{$1$} (G);
\path (E) edge[pil,color=red] node[right,black]{$2$} (H);
\path (F) edge[pil,color=blue] node[ right,black]{$1$} (I);

\path (G) edge[pil,color=red] node[left,black]{$2$} (J);
\path (H) edge[pil,color=blue] node[ right,black]{$1$} (K);
\path (I) edge[pil,color=red] node[below right,black]{$2$} (K);
\path (I) edge[pil,color=blue] node[ right,black]{$1$} (L);

\path (J) edge[pil,color=red] node[left,black]{$2$} (M);
\path (K) edge[pil,color=blue] node[ right,black]{$1$} (N);
\path (L) edge[pil,color=red] node[below right,black]{$2$} (N);

\path (M) edge[pil,color=blue] node[below left,black]{$1$} (O);
\path (N) edge[pil,color=red] node[right,black]{$2$} (O);
\end{tikzpicture}

\hspace{20mm}
\begin{tikzpicture}[>=latex]
\node (A) {\ytableaushort{\none {1}{1},{1}{1}}};
\node[below= 1  of A] (B) 
{\ytableaushort{\none {1}{2},{1}{1}}};
\node[below= 1  of B] (C) 
{\ytableaushort{\none {2}{2},{1}{2}}};
\node[below right = 1 and 0.5 of B] (D) 
{\ytableaushort{\none {1}{3},{1}{1}}};

\node[below= 1  of C] (E) 
{\ytableaushort{\none {2}{2},{2}{2}}};
\node[below= 1  of D] (F) 
{\ytableaushort{\none {2}{3},{1}{2}}};

\node[below= 1  of E] (G) 
{\ytableaushort{\none {2}{3},{2}{2}}};
\node[below= 1  of F] (H) 
{\ytableaushort{\none {3}{3},{1}{3}}};

\node[below= 1  of G] (I) 
{\ytableaushort{\none {3}{3},{2}{3}}};
\node[below= 1  of I] (J) 
{\ytableaushort{\none {3}{3},{3}{3}}};

\path (A) edge[pil,color=blue] node[ left,black]{$1$} (B);
\path (B) edge[pil,color=blue] node[ left,black]{$1$} (C);
\path (B) edge[pil,color=red] node[above right,black]{$2$} (D);

\path (C) edge[pil,color=blue] node[ left,black]{$1$} (E);
\path (C) edge[pil,color=red] node[above right,black]{$2$} (F);
\path (D) edge[pil,color=blue] node[right,black]{$1$} (F);

\path (E) edge[pil,color=red] node[left,black]{$2$} (G);
\path (F) edge[pil,color=blue] node[above left,black]{$1$} (G);
\path (F) edge[pil,color=red] node[right,black]{$2$} (H);

\path (G) edge[pil,color=red] node[left,black]{$2$} (I);
\path (H) edge[pil,color=blue] node[above left,black]{$1$} (I);
\path (I) edge[pil,color=red] node[left,black]{$2$} (J);
\end{tikzpicture}

}
\]
\caption{Two of the  connected components in the crystal graph of $\sr^3(\lm)$ with $\l=(3,2)$ and $\m=(1)$.}\label{dnk-2-1}
\end{figure}

 \subsection{Highest weight elements and  lattice words}

In this subsection, we use lattice words to give  a characterization of when a SVRPP is highest weight.
Let us first defined the  {\it reading word} of a SVRPP $T$, denoted   $\mathrm{read}(T)$. For the sake of description, we call the minimum number in each box of $T$ a {\it pivot}.   Then  $\mathrm{read}(T)$ is
the sequence  by reading the numbers of $T$ row by row, bottom to top, according to the following rule within each row: 
\begin{itemize}
    \item First, ignore all  the pivots in the row, and read the remaining numbers 
   from right to left and from largest to smallest within each
 box.  Then, read the pivots in this row from left to right, but when a pivot  is equal to the largest number  in the box  directly above, it will be ignored.  
\end{itemize}
See Figure \ref{Re-09} for an example, where the ignored pivots are colored red, and we use vertical lines to divide  subwords  from different rows, and in each subword we use a  dot to  divide   pivot and non-pivot  elements.
 \begin{figure}[h t]
\begin{center}
\begin{tabular}{c}
\begin{ytableau}
\none& 12 & 3& 4 & 567\\
123 & 3 &\ml{3}4 &\ml{4}\\
\ml{3} & \ml{3} & \ml{4}5 & 67\\
\ml{3}4 & 56
\end{ytableau}  
\end{tabular}
\end{center}
\caption{$\mathrm{read}(T)=64\cdot 5\,|\, 75\cdot 6\,|\, 432\cdot 13\,|\, 762\cdot 1345$.}\label{Re-09}
\end{figure}

\begin{remark}
The reading word of a SVRPP incorporates the constructions of  the reading word of a set-valued tableau and the reading word of a reverse plane partition defined   in \cite{Bandlow}. When restricting to   reverse plane partitions, our reading word is different from that in  \cite{Galashin 2}, where the reading word is defined by ignoring every  number equal to the number directly below it. 
\end{remark}

A word $w = w_1w_2 \cdots w_r$   is called a {\it reverse
lattice word}  (or a {\it Yamanouchi word}) if  each (backward) initial subword $w_k\cdots w_r$ has at least as many $a$’s as
$b$’s whenever $a<b$. For example, $312211$ is a reverse lattice word, while   $31221$ is not.

\begin{proposition}\label{GUk-09}
Let  $T\in \sr^n(\lm)$. Then $T$   is a highest weight element if and only if $\mathrm{read}(T)$ is a reverse lattice word.
\end{proposition}

To prove Proposition \ref{GUk-09}, we first list the following lemma. 

\begin{lemma}\label{njk-9-3-6}
Let  $T\in \sr^n(\lm)$, and for $1\leq i<n$, let    $T^{(i)}$ be the subtableau of $T$ occupied by $i$ and $i+1$. Then the subsequence of $\mathrm{read}(T)$ formed by $i$ and $i+1$ is equal to $\mathrm{read}(T^{(i)})$.
\end{lemma}

\begin{myproof}
Suppose that for $j\geq 1$, $u^{(j)}$   (resp., $v^{(j)}$) is the subword  in $\mathrm{read}(T)$  (resp., $\mathrm{read}(T^{(i)})$)  contributed by $i$'s and $i+1$'s in the $j$-th row of $T$  (resp., $T^{(i)}$). We need to verify  that $u^{(j)}=v^{(j)}$. This can be done by checking the following cases separately.  
\begin{itemize}
    \item[(1)] All $i$ and $i+1$ in the $j$-th row of $T$ are pivots. 

    \item[(2)] Otherwise, there exists at least one $i$ or $i+1$ in the $j$-th row of $T$ that is not a pivot. Precisely, we have the following subcases. 

    \begin{itemize}
        \item[(2.1)] Only the leftmost $i$  in the $j$-th row of $T$ is not a pivot. 

\item[(2.2)] All $i$'s are   pivots, and some (unique) $i+1$  in the $j$-th row of $T$ is not a pivot.

        \item[(2.3)] Both the leftmost $i$ and some (unique) $i+1$  in the $j$-th row of $T$ are not  pivots.
    \end{itemize}
\end{itemize}

According to the reading rule, it can be directly checked that in each of the above situations, there holds that 
$u^{(j)}=v^{(j)}$.  For example, we look at the case in (2.3). The distributions of $i$ and $i+1$ in the $j$-th row of $T$ are illustrated in Figure \ref{Rederu-09-2}, where the $*$'s stand for numbers not equal to $i$ and $i+1$. The colored numbers are not pivots. Note that  the leftmost non-pivot $i$ in $T$ becomes a  pivot in $T^{(i)}$.  By the reading rule,  we see that $u^{(j)}=\color{blue}i+1\, i\cdot w$ and $v^{(j)}=\color{blue}i+1\cdot i \, w$, where we use a dot to divide pivot and non-pivot elements, and the sequence $w$ is formed by the remaining pivots that are not ignored.  
\begin{figure}[h t]
\begin{center}
\renewcommand{\arraystretch}{1.5}
\begin{tabular}{|c|c|c|c|c|c|}
\hline
    $\ast\ast {\color{blue}i}$ & $i$ & $i$ & $i$,\  ${\color{blue}i+1} $ &$i+1 $& ${i+1} \ast\ast$\\
\hline 
\end{tabular}    
\end{center}
\caption{Illustration of the distribution of $i$ and $i+1$ in (2.3).}\label{Rederu-09-2}
\end{figure}
 \end{myproof}

 We now give a proof of Proposition \ref{GUk-09}.

\noindent
{\it Proof of Proposition \ref{GUk-09}.}
Assume that $T$ is a highest weight element, that is, $e_i(T)=0$ for every $i\in [n-1]$. This means that all $(i+1)$-pure columns of $T$ are paired with some $i$-pure columns. 
Consider the subsequence  
in $\mathrm{read}(T)$ consisting  of $i$ and $i+1$. By Lemma \ref{njk-9-3-6}, this subsequence coincides with $\mathrm{read}(T^{(i)})$. We explain that  each $i+1$ in $\mathrm{read}(T^{(i)})$ is paired with an $i$  located to its right, and different $i+1$ is paired with different $i$.  
\begin{itemize}
    \item If the $i+1$ is contributed by a mixed column of $T^{(i)}$, then we select the (topmost) $i$ contributed by this mixed column. 

    \item If the $i+1$ is contributed by an $(i+1)$-pure column of $T^{(i)}$, then we select the (topmost) $i$ contributed by the unique $i$-pure column that pairs with this  $(i+1)$-pure column in the $i$-signature of $T^{(i)}$.
\end{itemize}
Therefore, in each (backward) initial subword of $\mathrm{read}(T^{(i)})$, the number of $i$'s is always greater than or equal to the number of $i+1$'s. This in turn verifies  that  $\mathrm{read}(T)$ is a reverse lattice word.

We next prove the converse direction. Let   $\mathrm{read}(T)$ be a reverse lattice word. Suppose to the contrary that $T$ is not a highest weight element. Then there exists an $ i\in [n-1]$ such that $e_i(T)\neq 0$.
Consider the subtableau $T^{(i)}$.
 Let $A_i$  be the $(i+1)$-pure column corresponding to the leftmost unpaired  $\bminus$. We consider 
 the position of   the $i+1$ in $\mathrm{read}(T^{(i)})$ contributed by $A_i$.
 We have the following observations. 
 \begin{itemize}
     \item The  
$i$ contributed by a column (no matter it is $i$-pure or mixed) to the left of $A_i$  must be located  to the left of this $i+1$.

\item The 
$i$ (resp., $i+1$) contributed by an $i$-pure (resp., $(i+1)$-pure) column  to the right of $A_i$  must be located  to the right of this $i+1$. Moreover, by the choice of $A_i$, each such $i$ must be paired with some $i+1$ contributed by an $(i+1)$-pure column to the right of $A_i$.

\item   Both the  
$i$ and $i+1$ contributed by a mixed column  to the right of $A_i$  must be located  to the right of this $i+1$.
\end{itemize}
By the above observations, we see that 
 the (backward) initial   subword of $\mathrm{read}(T^{(i)})$, which    starts with the $i+1$ contributed by $A_i$, has at least one more $i+1$ than $i$.  In view of Lemma \ref{njk-9-3-6}, we conclude that  $\mathrm{read}(T)$ is not a reverse lattice word, leading to a   contradiction.
\qed

\begin{remark}\label{rer-9-23}
T. Yu  \cite{Yu-1}  suggested that one could also consider the column reading word $\mathrm{read}^c(T)$ defined as follows.
Read the columns from left to right, and within each column, remove the ingored pivots and then read the remaining numbers   in decreasing order. For example, the column reading word of $T$ in Figure \ref{Re-09} is $4321\ 65321\ 543\ 764\ 765$. 
Using similar arguments as given in the proof of Proposition \ref{GUk-09}, we can show that 
 $T\in \sr^n(\lm)$ is a highest weight element if and only if $\mathrm{read}^c(T)$ is a reverse lattice word.

\end{remark}

\subsection{A reconstruction algorithm}

In this subsection, we show  that  a SVRPP $T$ can be uniquely determined by  its reading word $\mathrm{read}(T)$,   its excess statistic $\mathrm{ex}(T)$, as well as its  \textit{height vector} $\mathrm{height}(T)$. Here the  {height vector} $\mathrm{height}(T)$ of $T$ is the sequence recording  the row indices  of the entries in  $\mathrm{read}(T)$. In other words, if write $\mathrm{read}(T)=w_1 w_2\cdots w_m$ and $\mathrm{height}(T)=(h_1,h_2,\ldots, h_m)$, then $w_i$ lies in row $h_i$ for $1\leq i\leq m$. For example, for the SVRPP $T$  in Figure \ref{Re-09}, we have $$\mathrm{height}(T)=(4,4,4,3,3,3,2,2,2,2,2,1,1,1,1,1,1,1).$$

\begin{remark}
We can similarly define the reading word and height vector of a benign tableau. Note that the descent-resolution operation   preserves the height vector. Thus the operators $e_i$ and $f_i$ on $\sr^n(\lm)$ also preserve  the height vector.   
\end{remark}

\begin{proposition}
 \label{lemma:injective}
 Fix a skew shape $\lm$, and  sequences $w, h$ and $\alpha$. Then   there exists at most one SVRPP $T$ of shape $\lm$ with $\mathrm{read}(T)=w$, $\mathrm{height}(T)=h$ and $\mathrm{ex}(T)=\alpha$.
\end{proposition}

\begin{myproof}
Suppose that there exists a SVRPP $T$ of shape $\lm$ with $\mathrm{read}(T)=w$, $\mathrm{height}(T)=h$ and $\mathrm{ex}(T)=\alpha$. 

Assume that $\lambda/\mu$ has $N$ rows. Write $\alpha=(\alpha_1,\ldots, \alpha_N)$. 
According to the sequence $h$, we can divide   $w=w^{(N)}\cdots w^{(1)}$ into  subwords such that for $1\leq i\leq N$,  $w^{(i)}$  consists of the numbers contributed by  row $i$. Let us proceed to divide each subword $w^{(i)}$ into two parts,  according to the  excess value $\alpha_i$. Precisely, we write $w^{(i)}=u^{(i)}\cdot v^{(i)}$, where $u^{(i)}$  includes the initial $\alpha_i$ entries of $w^{(i)}$ (which are exactly the non-pivot elements in row $i$) and $v^{(i)}$ includes the pivots that are not ignored in row $i$. Note that   $u^{(i)}$ is strictly decreasing, and $v^{(i)}$ is weakly increasing. 
For the SVRPP $T$ in Figure \ref{Re-09}, we have
\[\mathrm{read}(T)=w=w^{(4)}w^{(3)}w^{(2)}w^{(1)}=64\cdot 5\,|\, 75\cdot 6\,|\, 432\cdot 13\,|\, 762\cdot 1345,\]
where, as before, we use vertical lines to divide subwords $w^{(i)}$, and in each subword we use a  dot to further divide  $u^{(i)}$ and $v^{(i)}$. 
 
We shall reconstruct  $T$ row by row from top to bottom.  
For $i=1,2,\ldots, N$, we recover the $i$-th row of $T$ as follows. 

First, we use $v^{(i)}$ to determine the pivot in each box of row $i$. If $v^{(i)}$ is empty, then fill each box   in row $i$  with the maximum number in the box directly above  it. If $v^{(i)}$ is not empty, then implement the following procedure. 
\begin{itemize}
\item[] Assume that $v^{(i)}=v_1v_2\cdots v_p$. 
    \item[1.] Set $j=p$; 

     \item[2.] Let $B$ be the rightmost box in row $i$ of $\lm$ which is not filled with any number. Assume that  $a$ is the maximum value in the box directly above it (if there is no such box, then  put $a:=0$);

     \item[3.] If $v_j>a$, then place   $v_j$ into the box $B$. If this happens, we decrease $j$ by 1, and go back to Step 2;

     \item[4.] Otherwise, place the number $a$ into $B$, and go back to Step 2.
\end{itemize}
When all the entries   $v_p,v_{p-1},\ldots, v_1$ have been exhausted, it is possible that there   exist some boxes in row $i$ that are still empty. For each of such boxes, we fill it with the number $a$ (here $a$ is as defined in Step 2). 

After the above procedure, we recover all the pivots in row $i$.
Next, we use $u^{(i)}$ to determine the ``extra" numbers in each box of row $i$. Suppose that there are $q$ boxes in row $i$, and  the pivots are $a_1\leq a_2\leq \cdots\leq a_q$ listed from left to right. 
For $1\leq j\leq q$, we put the elements in $u^{(i)}$, belonging to  the interval $(a_j, a_{j+1}]$ (here we set $a_{q+1}=\infty$), into the $j$-th box counted from left. 
\end{myproof}


\subsection{Schur expansion and  proof of Theorem \ref{THM-2}}

Theorem \ref{THM-2} will be a consequence of the  the following result. 

\begin{theorem}\label{LJF_09}
For any  highest weight element  $T \in \mathcal{B}=\sr^n(\lm)$, the connected component containing  $T$   is isomorphic to $\mathcal{B}(\nu)$, where $\nu = \mathrm{wt}(T)$.
Moreover, we have
\[
\mathrm{SVRPP}^n(\lm) \cong \bigoplus_{\nu } \mathcal{B}(\nu)^{\oplus H_{\lm,n}^{\nu}},
\]
where $H_{\lm,n}^{\nu}$ is the number of highest weight elements in $\HW(\mathcal{B})$ which have weight $\nu$, and so 
$${G}_{\lm}(\mathbf{x}_n;\mathbf{1}; \mathbf{1})=\sum_{T\in \HW(\mathcal{B})}s_{\mathrm{wt}(T)}(\mathbf{x}_n)=\sum_{\nu}H_{\lm,n}^{\nu}s_{\nu}(\mathbf{x}_n).$$
\end{theorem}

Since the crystal  operators $e_i$ and $f_i$
preserve both the $\mathrm{ceq}$ 
statistic  and the  $\mathrm{ex}$ statistic,  all SVRPP's in the same connected component  have equal  $\mathrm{ceq}$    and   $\mathrm{ex}$. 
Therefore, Theorem \ref{LJF_09}  admits the following refinement.

\begin{theorem}\label{CHI-07}
For any skew shape $\lambda/\mu$, 
$${G}_{\lm}(\mathbf{x}_n;\mathbf{t};\mathbf{w})=\sum_{\gamma}\sum_{\theta}\mathbf{t}^\gamma\mathbf{w}^\theta\sum_{\nu}H_{\lm,n}^{\nu,\gamma,\theta}\,s_{\nu}(\mathbf{x}_n),$$
where $\gamma$ and $\theta$ are over all weak compositions, $\nu $ is over all partitions, and  the coefficient $H_{\lm,n}^{\nu,\gamma,\theta}$ counts the number of highest weight elements $T\in \sr^n(\lm) $ of weight $\nu$ and with   $\ceq(T) = \gamma$ and $\ex(T)=\theta$.  
\end{theorem}

Combining Theorem \ref{CHI-07} and Proposition  \ref{GUk-09} leads us to a proof of Theorem \ref{THM-2}. 

\begin{theorem}[=Theorem \ref{THM-2}]
 The coefficient  $H_{\lm,n}^{\nu,\gamma,\theta}$      is equal to the number of  SVRPP's $T\in \sr^n(\lm) $  of weight $\nu$ and with   $\ceq(T) = \gamma$ and $\ex(T)=\theta$, such that the reading word of $T$  
is a  reverse lattice word.
\end{theorem}

\begin{remark}
 Letting $n$ tend to $\infty$, we have the limit version:
$${G}_{\lm}(\mathbf{x} ;\mathbf{t};\mathbf{w})=\sum_{\gamma}\sum_{\theta}\mathbf{t}^\gamma\mathbf{w}^\theta\sum_{\nu}H_{\lm}^{\nu,\gamma,\theta}s_{\nu}(\mathbf{x}),$$ 
where the coefficient $H_{\lm,n}^{\nu,\gamma,\theta}$ is the number of  SVRPP's $T\in \sr(\lm) $  of weight $\nu$ and with   $\ceq(T) = \gamma$ and $\ex(T)=\theta$, such that the reading word of $T$  
is a  reverse lattice word.   
\end{remark}

The rest   is devoted to a proof of Theorem \ref{THM-2}. To this end, we show that each connected component in the crytal graph of $\sr^n(\lm)$ is isomorphic to  $\mathcal{B}(\nu)$, where $\nu$ is the weight of the highest weight element in this component. The key ingredient  is to prove that the operators $E_i$ on   words and  the operators $e_i$ on SVRPP's   commute with the operation of taking reading word.

\begin{proposition}
 \label{lemma:intertw}
 Let $T$ be a SVRPP of shape $\lambda/\mu$. Then 
\begin{equation}\label{hf-0-86}
  E_i\left(\mathrm{read}(T)\right)=\mathrm{read}\left(e_i(T)\right).  
\end{equation}
 In particular,   $e_i(T)=0$ if and only if $E_i(\mathrm{read}(T))=0$.
\end{proposition}

\begin{myproof}
Recall that when applying   the operator $e_i$ to $T$, we locate the column of $T$ corresponding to  the leftmost unpaired  `$\bminus$' in the $i$-string of $T$,  and then replace every $i+1$ in this column by  $i$. This results in a new filling, denoted  $T'$. Then we  apply the descent-resolution algorithm
to the subfilling of  $T'$ occupied by $i$ and $i+1$, producing a SVRPP $T''
 =e_i(T)$. 
 
Clearly, we can define the reading word of $T'$ or $T''$ in the same way as a SVRPP. We confirm  \eqref{hf-0-86} by verifying the following   two statements:
\begin{itemize}
    \item[(I)] $\mathrm{read}(T')=\mathrm{read}(T'')$;
    \item[(II)] $E_i(\mathrm{read}(T))=\mathrm{read}(T')$.
\end{itemize}


Let us first check the    statement in (I).  
Assume that 
\[T'=P_0\rightarrow P_1\rightarrow \cdots \rightarrow P_m=T'',\ \ \ \text{where $m\geq 0$}.\] 

According to  the proof of    Proposition \ref{kou98}, each right arrow represents a type (M1) $\rightarrow$ (1M) descent resolution. 
It is not hard to check that  for each $0\leq t<m$, $\mathrm{read}(P_t)$ and $\mathrm{read}(P_{t+1})$ coincide.  This verifies  the    statement in (I).

It remains to prove the statement in (II). 
From the definition of the reading word, we observe  that  $\mathrm{read}(T')$ is obtained from $\mathrm{read}(T)$ by replacing an   $i+1$, which is contributed by the topmost $i+1$ in the $(i+1)$-pure column corresponding to the leftmost unpaired $\bminus$ in the $i$-signature of $T$, by an  $i$. Therefore, the second statement is equivalent to showing  that 
\begin{itemize}
    \item[(II')] 
an $(i+1)$-pure column of $T$ is labeled by the leftmost unpaired  $\bminus$   if and only if the corresponding entry $i+1$   in $\mathrm{read}(T)$  is also labeled by the leftmost unpaired  $\bminus$. 
\end{itemize}

By Lemma \ref{njk-9-3-6}, we may restrict $T$ to  a SVRPP  that are filled  with subsets of $\{i, i+1\}$.  
For an  $i$-pure  or $(i + 1)$-pure  column $A$ of $T$, 
the topmost $i$ or $i+1$ will contribute to $\mathrm{read}(T)$, and let 
$R(A)$ be   this $i$ or $i+1$ appearing in  $\mathrm{read}(T)$. When  $A$ is mixed,  both the topmost $i$ and the topmost $i+1$ will contribute to $\mathrm{read}(T)$, and  we use  $R^{(i)}(A)$ (resp., $R^{(i+1)}(A)$) to be  this  $i$ (resp., $i+1$) appearing in $\mathrm{read}(T)$. 

We next to check the following  three facts.
\begin{itemize}
    \item[(i)] If  the $\bminus$ for an $(i+1)$-pure column $A$ is unpaired, then the entry $R(A)$ is unpaired in $\mathrm{read}(T)$.

        \item[(ii)] 
If  the $\bminus$ for an $(i+1)$-pure column $A$ is paired to the $\bplus$ for some $i$-pure column $B$, then   $R(A)$ is also
paired  to some $i$ (not necessarily to $R(B)$) in $\mathrm{read}(T)$.
    
    \item[(iii)]If a column $A$ is mixed, then   $R^{(i+1)}
(A)$ is paired to
some $i$ (not necessarily to $R
^{(i)}(A))$ in $\mathrm{read}(T)$.

\end{itemize}

The above facts tell that the  $(i+1)$-pure columns of $T$ corresponding to unpaired  $\bminus$'s   exactly give rise to the unpaired  $i+1$'s  in $\mathrm{read}(T)$. Moreover, notice  that  the reading
word rule preserves the order of pure columns from left to right. So the statement in (II') will be true once the facts in (i), (ii) and (iii) are proved. 

Before we proceed, the following  observation is needed. 
\begin{itemize}
    \item[(O)] Let $A$ be an $(i+1)$-pure column of $T$, and $B\neq A$ be any other column of $T$.  Then,   if $B$
lies  to the left (resp., right) of $A$, then the   $R(B)$ (or, $R^{(i)}(B)$ or $R^
{(i+1)}(B)$ if $B$ is mixed)
appears also to the left (resp., right) of $R(A)$.
\end{itemize}

This observation can be readily checked  if $B$ has no box filled with $\{i, i+1\}$. Consider now that $B$ has a (unique) box filled with $\{i,i+1\}$. If $B$ lies to the left of $A$, it is easy to see from the reading rule that both $R^{(i)}(B)$ or $R^
{(i+1)}(B)$ appear to the left of $R(A)$. If $B$ lies to the right of $A$, the box filled with $\{i, i+1\}$  must lie in a row strictly above the topmost $i+1$ in  $A$, and so  both $R^{(i)}(B)$ or $R^
{(i+1)}(B)$ also appear to the right of $R(A)$.

We now give proofs of (i), (ii) and (iii). 
\begin{itemize}

    \item[(i)]  In this case, the $\bminus$ for $A$ is not paired.    Now, the $\bplus$ for every $i$-pure column
to the right of $A$ is paired to the $\bminus$ for some $(i+1)$-pure column. So, by the observation (O), for  every $i$ to the right of $R(A)$ in  $\mathrm{read}(T)$, there exists an $i+1$
between $R(A)$ and this $i$, and moreover for different $i $’s, these $i+ 1$’s are also
different. This leads to an injective map that attaches to each
$i$ lying to the right of $R(A)$  an $i+1$ between  it and $R(A)$. Therefore, in $\mathrm{read}(T)$, every $i$ to the right of $R(A)$ will be  paired to some $i+1$ lying to the right of $R(A)$. This implies that $R(A)$ must be  unpaired in $\mathrm{read}(T)$.

\item[(ii)] In this case, the $\bminus$ for $A$ is paired with the $\bplus$ for some $i $-pure column $B$ to the right of $A$. Let $C$ be the leftmost
$(i+1)$-pure column to the right of $B$ . Let us consider the open segment, denoted $(R(A), R(C))$, of $\mathrm{read}(T)$ between $R(A)$ and $R(C)$ (not including $R(A)$ and $R(C)$) (If there is no such $C$, we consider the whole segment to the right of $R(A)$). By (O), the elements in this segments are exactly contributed by the columns between $A$ and $C$. 
By the choices, we have an equal number  of $i's$ and $i+1's$ in $(R(A), R(B))$, while  in $[R(B),R(C))$, the number of $i's$ is strictly larger than the number of $i+1's$. So $R(A)$ must be paired with some $i$ in $(R(A), R(C))$.

\item[(iii)]  In this case, $A$ is mixed.  Let $B$ be the leftmost $(i+1)$-pure column to the
right of $A$. Consider the open segment $(R^{(i+1)}(A), R(B))$ in $\mathrm{read}(T)$ (If there is no such $B$, we consider the whole segment to the right of $R^{(i+1)}(A)$). Note that $R^{(i)}(A)$ belongs to $(R^{(i+1)}(A), R(B))$. 

Let $C$ be any column other than $A$ and $B$. We assert that if $C$ contributes an $i+1$ to $(R^{(i+1)}(A), R(B))$, then it must be mixed and moreover it   contributes an $i$ to $(R^{(i+1)}(A), R(B))$. This can be seen as follows.   By the observation in (O) and the choice of $B$, it is easy to confirm  that $C$ cannot be $(i+1)$-pure. So $C$ is mixed. By our reading rule, it can be  easily checked that $C$ must lie between $A$ and $B$. In this case, again by our reading rule, we can check that $R^{(i)}(C)$ must belong  to $(R^{(i+1)}(A), R(B))$. This concludes the assertion. 

Now we see that in   $(R^{(i+1)}(A), R(B))$, the number of $i$'s is strictly larger than the number of $i+1$'s (since $A$ contributes an extra $i$). This means that $R^{(i+1)}(A)$ must be paired with some $i$ in the segment $(R^{(i+1)}(A), R(B))$. 
\end{itemize}
This completes the proof of  the statement in (II').
\end{myproof}

Proposition \ref{lemma:intertw} leads to the following commuting relation for $F_i$ and $f_i$.

\begin{corollary}\label{ago-9-0}
 Let $T$ be a SVRPP of shape $\lambda/\mu$. Then 
\begin{equation}\label{dhudi-68}
 F_i\left(\mathrm{read}(T)\right)=\mathrm{read}\left(f_i(T)\right).   
\end{equation} 
Here, we also have   $f_i(T)=0$ if and only if $F_i(\mathrm{read}(T))=0$.
\end{corollary}

\begin{myproof}
When $f_i(T)\neq 0$,  \eqref{dhudi-68}   follows from \eqref{hf-0-86} since both $e_i/f_i$ and $E_i/F_i$ are crystal operators. We next explain that $f_i(T)=0$ if and only if $F_i(\mathrm{read}(T))=0$.  In the  proof of Proposition \ref{lemma:intertw}, we have seen that the  $(i+1)$-pure columns of $T$ corresponding to unpaired  $\bminus$'s   exactly give rise to the unpaired  $i+1$'s  in $\mathrm{read}(T)$, that is, $\varepsilon_i(\mathrm{read}(T))=\varepsilon_i(T)$. Notice that $S_{n}^m$ has a seminormal crystal structure and $\mathrm{wt}(T)=\mathrm{wt}(\mathrm{read}(T))$.
From these facts and Corollary \ref{gsgl-3}, we have  $$\varphi_i(\mathrm{read}(T)) - \varepsilon_i(\mathrm{read}(T))=\mathrm{wt}(T)_i-\mathrm{wt}(T)_{i+1}=\varphi_i(T) - \varepsilon_i(T),$$
from which we deduce that $f_i(T)=0$ if and only if $F_i(\mathrm{read}(T))=0$.
\end{myproof}

We are finally in a position to give a proof of Theorem \ref{LJF_09}.

\noindent
{\it Proof of Theorem \ref{LJF_09}.} 
Fix a skew shape $\lm$,  a   sequence   $h$ of positive integers,  and a composition $\theta$ of nonnegative integers. Suppose that $h$ has $m$ entries. By  Proposition  \ref{lemma:injective} and Corollary \ref{ago-9-0}, we get   an injection  $T\mapsto \mathrm{read}(T)$ from the set of  SVRPP's $T\in \sr^n(\lm) $  with $h(T)=h$ and $\mathrm{ex}(T)=\theta$ to  the set $S_{n}^m$. This map  takes the operators $e_i$ and $f_i$ to $E_i$ and $F_i$, respectively. 

Now, each connected component in  the crystal  graph of $\sr^n(\lm) $ is 
embedded into a connected component in the crystal graph of $S_{n}^m$. Moreover, because $e_i(T)=0 \Leftrightarrow E_i(\mathrm{read}(T))=0$ and $f_i(T)=0\Leftrightarrow F_i(\mathrm{read}(T))=0$, we easily deduce that the embedding is actually surjective. This, together with Proposition \ref{lemma:crystal}, completes the proof of Theorem \ref{LJF_09}. 
\qed

\section{Proof of Theorem \ref{Newton-98}}\label{Sec-5--2}

We begin with a criterion of when a symmetric polynomial has saturated Newton polytope  given in \cite{Nguyen}.
As usual, for two partitions $\lambda$ and $\lambda'$, write $\lambda \supseteq \lambda'$ if $\lambda_i \geq \lambda_i'$ for every $i$.

\begin{theorem}[\cite{Nguyen}]\label{JU-0-j}
Let $F(\mathbf{x}_n)$ be a linear combination of Schur polynomials.  If
\begin{itemize}
\item[(1)] the Schur polynomials in each  homogeneous component  of  $F(\mathbf{x}_n)$ have the same sign;

    \item[(2)] there exist partitions $\mu, \nu$ such that $s_\rho(\mathbf{x}_n)$ appears in $F(\mathbf{x}_n)$ if and only if $\mu \subseteq \rho\subseteq \nu$,
\end{itemize}
then  $F(\mathbf{x}_n)$ has SNP.
\end{theorem}

\begin{remark}\label{ahhn-6}
    Let $\mathcal{P}$ be a lattice polytope, namely, the convex hull of some lattice points in $\mathbb{Z}^n$. We say that $\mathcal{P}$ has \textit{integer decomposition property}  (IDP) if, for any positive integer $t$ and any lattice point $p$ in the dilation $t\mathcal{P}=\{t v \colon v \in \mathcal{P}\}$, there are $t$ lattice  points $v_1, \dots, v_t \in \mathcal{P}$ such that $p=v_1+\dots+v_t$. 
    To determine whether a lattice polytope has IDP is a basic question in   Ehrhart theory.    
    It was shown in \cite{Nguyen} that if  $F(\mathbf{x}_n)$ satisfies the conditions in Theorem \ref{JU-0-j}, then  its Newton polytope $\mathrm{Newton}(F(\mathbf{x}_n))$ has IDP.
\end{remark}

To give a proof of Theorem \ref{Newton-98}, we shall show that hybrid Grothendieck polynomials indeed satisfy the conditions in  Theorem \ref{JU-0-j}. Our next task is to determine  the partitions $\mu$ and $\nu$  needed in  Theorem \ref{JU-0-j}. The partition $\mu$ is simply defined as $\mu=(\lambda_1)$. 
We next construct the partition $\nu$.

Fix  a partition $\l=(\l_1,\l_2,\ldots)$. Let $\ell$ be the length of $ \lambda$, namely, the number of positive  parts  of $\lambda$. Define
$$\l^\flat=(\l^\flat_1,\l^\flat_2,\ldots)$$  
by setting 
$\l^\flat_1=\l_1$, and for $i\geq 2$,
\[
\l^\flat_i=\max\left\{0,\ \min\{\l^\flat_{i-1}-1,\l_i\}\right\}.
\]
It is easy to see that $\l^\flat$ is a partition contained in $ \l$. 
Moreover, the positive parts of $\l^\flat$ strictly decrease, that is, $\l^\flat_1>\l^\flat_2>\cdots>\l^\flat_k$, where $k$ is the length of $\l^\flat$. 
For example, for $\l=(7,5,4,4,1,1)$, we have $\l^\flat=(7,5,4,3,1)$ and so $k=5$. 

\begin{remark}\label{hgk-9-6}
We say that a partition is strict if its positive parts are strictly decreasing.  
It is readily verified  that $\l^\flat$ is the largest among all strict partitions contained in $\lambda$. In other words, for any strict partition $\rho\subseteq \lambda$, we have $\rho \subseteq \l^\flat$.
\end{remark}

For a fixed  positive integer $n$, define 
$$\bar{\l}^{(n)}=(\bar{\l}^{(n)}_1,\bar{\l}^{(n)}_2,\ldots),$$
where $\bar{\l}^{(n)}_i=0$ for $i>n$, and 
\begin{itemize}
    \item[$\bullet$] when $n\leq k$, set $\bar{\l}^{(n)}_i=\l^\flat_i+i-1$ for $i\leq n$;

    \item[$\bullet$] when $n> k$, set $\bar{\l}^{(n)}_i=\l^\flat_i+i-1$ for $i\leq k$, and $\bar{\l}^{(n)}_i=k$ for $k<i\leq n$.
\end{itemize}
Since $\l^\flat$ is strict, it is easily seen that $\bar{\l}^{(n)}$ is a partition of length $n$.
Still, for the partition $\l=(7,5,4,4,1,1)$, we see that  $\bar{\l}^{(3)}=(7,6,6)$
 and $\bar{\l}^{(7)}=(7,6,6,6,5,5,5)$.  As will be shown later, the partition $\bar{\l}^{(n)}$ is the desired partition $\nu$.

 We proceed to construct the  highest weight elements, denoted as $T_{\min}$ and $T_{\max}$, in $\sr^n(\l)$ whose weights reach $\mu=(\lambda_1)$ and $\nu=\bar{\l}^{(n)}$, respectively.  Let $T_{\min}$ be simply defined by setting all entries  equal to  $\{1\}$. Obviously, its reading word is the string with exactly $\lambda_1$ $1$'s, and thus is a reverse lattice word. So $T_{\min}$   is   a highest weight element with  $\mathrm{wt}(T_{\min})=(\l_1)$.

The construction of  $T_{\max}$ is via the following rule: for $(i,j)\in \lambda$,
$$T_{\max}(i,j) = 
\begin{cases}
    \{i\}, & \text{if } 1\leq i\leq n \text{ and } j<\l^\flat_i, \\
    [i,n]=\{i,i+1,\ldots,n\}, & \text{if } 1\leq i\leq n \text{ and } j=\l^\flat_i, \\
    \{n\}, & \text{otherwise} .
\end{cases}
$$
It is easily seen  that $T_{\max}$ is a SVRPP. For  $\l=(7,5,4,4,1,1)$,   Figure \ref{fig:my_label-089} gives $T_{\max}$ for $n=3$  and $n=7$, respectively. 
\begin{figure}[H]
    \centering
    \[
\ytableausetup{boxsize=2.4em}
\scalebox{0.8}{
    \begin{tabular}{ccc}
    \begin{ytableau}
    1 & 1 & 1 & 1 & 1 & 1 & [1,3]\\
    2 & 2 & 2 & 2 & [2,3]\\
    3 & 3 & 3 & [3,3]\\
    3 & 3 & 3 & 3\\
    3\\
    3
    \end{ytableau}   
         &\quad\quad  &
    \begin{ytableau}
    1 & 1 & 1 & 1 & 1 & 1 & [1,7]\\
    2 & 2 & 2 & 2 & [2,7]\\
    3 & 3 & 3 & [3,7]\\
    4 & 4 & [4,7]& 7 \\
    [5,7]\\
    7\\
    
     \end{ytableau} 
    \end{tabular}
    }
    \]
    \caption{ $T_{\max}$ for $\l=(7,5,4,4,1,1)$ and  $n=3, 7$, respectively.}
    \label{fig:my_label-089}
\end{figure}

\begin{lemma}\label{nknb-0-p}
For any partition $\lambda$ and positive integer $n$, $T_{\max}$ is a highest weight element in $\sr^n(\l)$ with $\mathrm{wt}(T_{\max})=\bar{\l}^{(n)}$.   
\end{lemma}

\begin{myproof}
We first check that $T_{\max}$ has weight equal to $\bar{\l}^{(n)}$.
Write $\mathrm{wt}(T_{\max})=(c_1,\ldots, c_n)$. 
We distinguish into two cases.
Let  $k$ denote  the length  of   $\l^\flat$.   

\vspace{0.5em}
Case 1: $n\leq k$. Consider $c_i$ for $1\leq i\leq n$. When $1\leq i<n$, by the construction of $T_{\max}$, we see that  each of the first $i-1$ rows has  one box containing $i$, and in row $i$, each of the leftmost  $\l^\flat_i$ boxes contains $i$. Moreover, these  boxes live  in different columns. This yields that $c_i=\l^\flat_i+i-1=\bar{\l}^{(n)}_i$.

    We next consider the situation  for $i=n$. With the same observation as above  for $1\leq i<n$, we already have  $\bar{\l}^{(n)}_n$ boxes, lying in distinct columns, that contain $n$. However, each box  in the skew shape $\lambda/(\l^\flat_1,\ldots, \l^\flat_n)$ is also filled with $\{n\}$.  We conclude the proof by explaining that for any box $B$  in $\lambda/(\l^\flat_1,\ldots, \l^\flat_n)$, there exists a box  $B'\in \lambda$, lying above $B$ in the same column, that  contains  the number  $n$. 
    This is clear when $B$ lies strictly below row $n$. 
    
    Now we assume that $B$ lies in  or above row $n$. Denote   $B=(i_0,j_0)$.  Since $\l^\flat_1=\lambda_1$, we have $i_0>1$. We explained in Remark \ref{hgk-9-6}  that $\l^\flat$ is the maximum strict partition contained in $\lambda$.  This forces that $\l^\flat_{i_0}=\l^\flat_{i_0-1}-1$.  We claim that there exists   $1\leq i'<i_0$ such that $\l^\flat_{i'}\leq j_0$. Suppose to the contrary that there is no such $i'$, then we have $\l^\flat_{i}> j_0$ for every $1\leq i<i_0$. However, this would yield that  $\l^\flat_{i_0}\geq j_0$, contrary to the fact that $\l^\flat_{i_0}<j_0$. This verifies the claim. So any $n$ in the box $B\in \lambda/(\l^\flat_1,\ldots, \l^\flat_n)$ has no contribution to the weight. Combining the above analysis, we conclude that $c_n=\bar{\l}^{(n)}_n$. 
 
 \vspace{0.5em}
 Case 2: $n\geq k$. When $1\leq i\leq k$, it is easily checked  that $c_i=\bar{\l}^{(n)}_i$. When $k<i<n$, the boxes that contain $i$ are exactly the rightmost boxes in the  rows of  $\l^\flat$, implying that $c_i=\bar{\l}^{(n)}_i=k$. When $i=n$, the argument is similar as that for $i=n$ in Case 1.

It remains to verify that  $T_{\max}$ is a highest weight element. By Remark \ref{rer-9-23},  it suffices to check that the  column reading word $\mathrm{read}^c(T_{\max})$  is a reverse lattice word. By the construction of $T_{\max}$, we see that the contribution of each column is of the form $m\cdots 21$ for some $m$.  For example, for the left picture in  Figure \ref{fig:my_label-089}, we have 
\[
 \mathrm{read}^c(T_{\max})=321\cdot   321\cdot 321\cdot 321\cdot 321\cdot 1 \cdot 321. 
\]
From this observation, it is obvious that   $\mathrm{read}^c(T_{\max})$  is a reverse lattice word.
\end{myproof}

\begin{lemma}\label{yi-0-p}
For any partition $\lambda$ and positive integer $n$, if a partition $\rho$ satisfies $(\l_1)\subseteq \rho \subseteq \bar{\l}^{(n)}$, 
then there exists a highest weight element in $\sr^n(\l)$ whose weight is equal to  $\rho$.
\end{lemma}

\begin{myproof}
If $(\l_1)\subseteq \rho \subseteq \bar{\l}^{(n)}$, then we have $\rho_1=\l_1$ and $\bar{\l}^{(n)}_i-\rho_i\geq 0$ for $2\leq i \leq n$. Write $d_i=\bar{\l}^{(n)}_i-\rho_i$.

We define an operator $L_i$ acting on $T\in \sr^n(\l)$ as follows. Choose the leftmost column containing $i$ of $T$, and delete all $i$ in this column. This possibly results in  empty boxes in this column. If this happens,  we fill each empty box with $\{i-1\}$. After the action of $L_i$, we get a new   $L_i(T)\in \sr^n(\l)$. 

Consider $$T^*=L_2^{d_2}\cdots L_n^{d_n}(T_{\max}),$$
where we set  $L_i^0=\mathrm{id}$. In the proof of Lemma \ref{nknb-0-p}, we have observed that each column of $T_{\max}$ contributes to the column reading word a string of the form $m\cdots 2 1$. In other words, the numbers appearing in each column of $T_{\max}$ form a set $[m]$.  From this observation, it is easy to obtain that 
\[
 \mathrm{wt}(T^*)=  \bar{\l}^{(n)}-(0,d_2,\ldots, d_n)=\rho. 
\]
Moreover, it is not hard to check that the column reading word of  $T^*$ is a reverse lattice word, and so $T^*$ is the required  highest weight element.
\end{myproof}

\begin{lemma}\label{yi-0-p-2}
For any partition $\lambda$ and positive integer $n$,  if $T\in \sr^n(\l)$ is a highest weight element with $\mathrm{wt}(T)=\rho$, then $(\l_1)\subseteq \rho\subseteq \bar{\l}^{(n)}$.  
\end{lemma}

\begin{myproof} 
First, notice that each box in the first row of $T$ must contain $1$ since otherwise   $\mathrm{read}(T)$ would not be a reverse lattice word. This implies that $\rho_1=\lambda_1$. The remaining is devoted to a proof of   $\rho\subseteq \bar{\l}^{(n)}$. We need the following claim.
\begin{itemize}
    \item Fix  $p\geq 1$. Then,  for any $q\geq p$, there are at most $p-1$ columns, each of which contains a $q$ lying  in the first $p-1$ rows of $T$.
\end{itemize}
The argument is by  induction on $p$.
The case for $p=1$ is trivial.  For  $p=2$, since each box in the first row must contain $1$, the only possible column containing $q$ in the first row is the last column of $T$. This verifies the claim for $p=2$.  

Assume  now that $p>2$ and the claim is true for  $p-1$. Suppose  to the contrary that  there are at least $p$ columns (say, columns $C_1,\ldots, C_p$) containing $q$ in the first $p-1$ rows of $T$. Consider column $C_j$ for each fixed $1\leq j\leq p$. 
\begin{itemize}
    \item[(1)] If $C_j$  contains (at least one) $q-1$, then the topmost $q-1$ either lies strictly above row $p-1$, or lies in  row $p-1$ and stays in the same box as a  $q$.  

    \item[(2)] If $C_j$  does not contain  $q-1$, then it provides a $\bminus $ in the $(q-1)$-signature. Since $T$ is a highest weight element,  there must exist  a column $C_j'$, which contains  $q-1$ (but not $q$) and lies  to the right of $C_j$, providing   a $\bplus$.  Moreover, the topmost $q-1$ in  column $C_j'$ lies strictly above row $p-1$. 
\end{itemize}
By the above analysis and noticing  that 
there exists at most one box   in row $p-1$ that contains  both $q-1$ and $q$,  we obtain that  there are at least $p-1$ columns containing $q-1$ lying  in first $p-2$ rows of $T$. This leads to   a contradiction, and so the claim is verified. 

Next, we   prove 
$\rho\subseteq \bar{\l}^{(n)}$ by contradiction. 
Suppose otherwise that  $\rho\nsubseteq \bar{\l}^{(n)}$.
Let  $i>1$ be the smallest index such that $\rho_i>\bar{\l}_i^{(n)}$. As before, let $k$ be the length of $\l^\flat$.  The discussion will be divided into two cases. 

Case 1: $i\leq \min\{k,n\}$.
 In this case, we have $\bar{\l}_i^{(n)}={\l}^{\flat}_i+i-1$. Keep in mind that $\rho_i$ is  the number of columns of $T$ that contain the number $i$. Because of $\rho_i>\bar{\l}_i^{(n)}$,  there are at least $i$ columns that contain $i$ and lie strictly to the right of column ${\l}^{\flat}_i$.
  Let $r=\max\{1\leq j\leq i-1\colon  {\l}^{\flat}_j-{\l}^{\flat}_{j+1}>1\}$ (set $r=0$ if there is no such $j$). 
 Obviously, we have ${\l}^{\flat}_{r+1}=\l_{r+1}$. So the columns, lying strictly to the right of column ${\l}^{\flat}_{r+1}$, are contained entirely in the first $r$ rows, and moreover it follows from  the choice  of $r$,  among these columns,  there are at least $r+1$  columns   containing $i$.  
Because $i\geq r+1$, this is contrary to the claim by setting $p=r+1$ and $q=i$.  
 

 Case 2: $k<i\leq n$. Let $\ell$ denote the length of $\lambda$. 
Note that $k\leq \ell$. 
We still distinguish into two situations. 
\begin{itemize}
    \item[(1)]  $k<\ell$. In this situation,  we must have ${\l}^{\flat}_k=1$ and so $\bar{\l}_k^{(n)}=k$. Now we  see that $\bar{\l}_i^{(n)}=\bar{\l}_k^{(n)}=k$. By Case 1, we have $\rho_k\leq\bar{\l}_k^{(n)}$. Together with $\rho_i\leq \rho_k$, we arrive at   $\rho_i\leq  \bar{\l}_i^{(n)}$, leading to a contradiction. 

    \item[(2)]  $k=\ell$. Clearly, one can restrict  $i=k+1$. However, the assuption  $\rho_{k+1}>\bar{\l}_{k+1}^{(n)}=k$ would contradict the claim by setting $p=q=k+1$.
\end{itemize}
This completes the proof that 
$\rho\subseteq \bar{\l}^{(n)}$. \end{myproof}


We can now finish the proof of 
Theorem \ref{Newton-98}.

\noindent
{\it Proof of Theorem \ref{Newton-98}.}
The proof follows from Theorem \ref{JU-0-j} together with 
Lemmas \ref{yi-0-p} and \ref{yi-0-p-2}.
\qed

\begin{remark}
    By Remark \ref{ahhn-6}, the Newton polytope of $G_{\l}(\mathbf{x}_n;\mathbf{1};\mathbf{1})$  has IDP.
\end{remark}

We do not know if a hybrid Grothendieck polynomial  $G_{\l/\mu}(\mathbf{x}_n;\mathbf{1};\mathbf{1})$ of skew shape has SNP. This seems unknown even for the stable or dual   stable Grothendieck polynomials of skew shapes.

As another consequence of   Lemmas \ref{yi-0-p} and \ref{yi-0-p-2}, we have the following corollary.

\begin{corollary}\label{bnf-5-8}
The degree of     $G_{\l}(\mathbf{x}_n;\mathbf{t};\mathbf{w})$
is equal to the size of $\bar{\l}^{(n)}$. 
\end{corollary}

\begin{remark}
The  degree of the stable Grothendieck polynomial $G_{\l}(\mathbf{x}_n;\mathbf{0};\mathbf{-1})$ has been considered in 
\cite{Rajchgot}. As a comparison, it should be noticed that  $G_{\l}(\mathbf{x}_n;\mathbf{0};\mathbf{-1})$ is nonzero if and only if $n$ is greater than or equal to the length of $\lambda$. However we do not have   such a constraint  for a hybrid Grothendieck polynomial, that is, $G_{\l}(\mathbf{x}_n;\mathbf{t};\mathbf{w})$ is nonzero for any $n\geq 1$.   
\end{remark}

\section{The omega involution image}\label{info-14}

Recall that 
\begin{align*}
{G}_{\lm}(\mathbf{x} ;\alpha;\beta)={G}_{\lm}(\mathbf{x} ;\mathbf{t};\mathbf{w})|_{t_i=\alpha, w_i=\beta}=\sum_{T\in \mathrm{SVRPP}(\lambda/\mu)} {\alpha}^{|\mathrm{ceq} (T)|}{\beta}^{|\mathrm{ex}(T)|}
\textbf{x}^{\mathrm{ircont}(T)},
\end{align*}
and ${J}_{\lm}(\mathbf{x} ;\alpha;\beta)=\omega \left({G}_{\lm}(\mathbf{x} ;\alpha;\beta)\right)$ is its image. 
Our goal in this section is to  give a combinatorial description of   ${J}_{\lm}(\mathbf{x} ;\alpha;\beta)$, thus proving Theorem \ref{THM-33}.  

Our starting point is to use  Fomin--Greene type operators to give a construction of   ${G}_{\lm}(\mathbf{x} ;\alpha;\beta)$. 
Fix a partition $\mu$. Let \[\bigoplus_{\mu\subseteq\lambda}\mathbb{Z}[\alpha, \beta] \{\lambda\} \] be the free $\mathbb{Z}[\alpha, \beta]$-module generated by the partitions containing $\mu$. 
We need two kinds  of operators acting on this   module. For $\mu\subseteq \lambda$, the {\it column adding operators} $u_i$  are defined by 
\[
 u_{i}\cdot\lambda=  \sum_{a\geq 1} \alpha^{a-1 } (\lambda \cup \left \{ \text{$a$ boxes in the $i$-th column of $\lambda$}   \right \}),  
\]
where  the sum is over all partitions by adding $a$ boxes to the $i$-th column of $\lambda$ (we set $u_{i}\cdot\lambda=0$ if there is no valid shape formed in this way). 
When $\alpha=1$, this kind of operators have appeared  in \cite{Lam} for constructing ordinary reverse plane partitions. 
To generate SVRPP's, we   borrow the {\it Schur operators}   \cite{Fom-2} with the parameter    $\beta$:
  \begin{align*}
         d_{\mu,i}\cdot \lambda=\begin{cases}
 \beta  (\lambda \setminus  \left \{\text{one box in the $i$-th column of $\lambda$}   \right \}), & \text{ if this gives a valid shape }  \\
 & \text{  containing $\mu$,}  \\[5pt]
 0, & \text{ otherwise. } 
\end{cases}
     \end{align*}
In Figure \ref{FO-098}, we give some examples  of the operators $u_i$ and $d_{\mu,i}$.
     \begin{figure}[h t]
     \begin{center}
\[
\ytableausetup{boxsize=1.5em}
\scalebox{0.9}{
\begin{tabular}{cccc}
& & & \\
\begin{tabular}{c}
\\
\\

$u_2 \enspace \cdot $
\end{tabular}

\ydiagram{3,1,1}
\begin{tabular}{c}
\\
\\

$=$
\end{tabular}

\ydiagram{3,2,1}
\begin{tabular}{c}
\\
\\

$+\ \alpha$
\end{tabular}
\ydiagram{3,2,2}\quad
\begin{tabular}{c}
\\
\\

,
\end{tabular}\\

& & & \\
& & & \\
\begin{tabular}{c}
\\
\\

$d_{(2,1,1),1} \enspace \cdot$
\end{tabular}
\ydiagram{2,1,1,1}
\begin{tabular}{c}
\\
\\

$=\ \beta$
\end{tabular}
\ydiagram{2,1,1}
\begin{tabular}{c}
\\
\\

$,\quad d_{(3,1,1),1} \enspace \cdot$
\end{tabular}

\ydiagram{3,1,1}
\begin{tabular}{c}
\\
\\

$=0$.
\end{tabular}
 \\
& & & \\
\end{tabular}
}
\]
\end{center}
\caption{Illustration  of operators $u_i$ and $d_{\mu,i}$.}\label{FO-098}
\end{figure}

Consider the operator
\[\tilde{u}_i=u_i(1+d_{\mu,i}).\]
Set
\[\tilde{A}(x)=\cdots (1+x\tilde{u}_3) (1+x\tilde{u}_2)(1+x\tilde{u}_1).\]

\begin{proposition}\label{YTRd-9}
The coefficient of   $\lambda$ in the product  \[\cdots \tilde{A}(x_3) \tilde{A}(x_2)\tilde{A}(x_1)\,\mu\] is equal to ${G}_{\lm}(\mathbf{x} ;\alpha;\beta)$.
\end{proposition}

\begin{myproof}
For $j\geq 1$, we have $\tilde{A}(x_j)= \cdots (1+x_j\tilde{u}_3) (1+x_j\tilde{u}_2)(1+x_j\tilde{u}_1)$. Let us see the effect of the  term $x_j\tilde{u}_i=x_ju_i+x_j u_i d_{\mu,i}$. 
\begin{itemize}
\item[(1)] $x_ju_i$:  we add any possible  vertical strip to the $i$-th column, and then fill   each box in the strip with the number $j$. If the added strip has $a$ boxes, then it will contribute a factor $\alpha^{a-1}x_j$. 

\item[(2)] $x_ju_id_{\mu,i}$: If the bottom box in the $i$-th column is removable, then we   remove this box, followed by  adding any possible vertical strip. Actually, this procedure has no effect to the removable box (it is removed, and then added back). We fill the removable box, as well as the remaining boxes in the added strip, with the number $j$.    If the added strip has $a$ boxes, then this procedure will contribute a factor $\alpha^{a-1}\beta x_j$. 
\end{itemize}
Notice that the operation in (2) will produce extra numbers in the removable boxes. The two procedures in (1) and (2) allow us to construct  SVRPP's. Conversely, each SVRPP can be produced uniquely via these two types of operations. Moreover, it is easily seen that each such constructed SVRPP has weight exactly equal to ${\alpha}^{|\mathrm{ceq} (T)|}{\beta}^{|\mathrm{ex}(T)|}
\textbf{x}^{\mathrm{ircont}(T)}$.
So the proof is complete. 
\end{myproof}

Our next task is to show that the operators $\tilde{u}_i$ satisfy conditions  required   in the   following classic result by    Fomin and Greene \cite{Fomin}.  

\begin{theorem}[\cite{Fomin}]\label{F}
Suppose that a set  $\left \{ v_i \colon i\in \mathbb{Z}   \right \} $ of elements in an associative
 algebra  satisfies the following relations:
    \begin{itemize}
        \item[(i)] $v_i v_k v_j=v_k v_i v_j$,\ \ \ \text{for $i<j<k$};
         \item[(ii)] $v_j v_k v_i=v_j v_i v_k$,\ \ \ \text{for $i<j<k$};
          \item[(iii)] $v_j(v_iv_j-v_jv_i)=(v_iv_j-v_jv_i)v_i$,\ \ \ \text{for $i<j$}.
    \end{itemize}
 Then there hold the following statements. 
\begin{itemize}
    \item[(1)] 
For $k\geq 0$, let 
\begin{align*}
e_k(\mathbf{{v}})=e_k(v_1,v_2,\ldots)=\sum_{a_1>a_2>\cdots >a_k}^{}  v_{a_1}v_{a_2}\cdots v_{a_k}
    \end{align*}
    be the noncommutative analogue of  elementary symmetric functions. Then $e_k(\mathbf{{v}})$ commute, that is, $e_i(\mathbf{{v}})\cdot e_j(\mathbf{{v}})=e_j(\mathbf{{v}})\cdot e_i(\mathbf{{v}})$ for any nonnegative integers $i$ and $j$. 

\item[(2)]
Let $A(x)$ and $B(x)$ be defined as 
    \begin{align*}
         A(x)=\cdots(1+xv_2)(1+xv_1),\ \ \ \  
    B(x)=\frac{1}{1-xv_1}\frac{1}{1-xv_2}\cdots. 
    \end{align*}
Assume that  $\mathbf{x}=(x_1,x_2,\ldots)$ is a sequence of   commuting indeterminates which 
commute  with each   $v_j$.  As before, write  $\lambda'$ for the conjugate of $\lambda$.     Then one has the following Cauchy type identities: 
\begin{align*}
         \cdots A(x_2)A(x_1)=\sum_{\lambda }  s_{\lambda}(\mathbf{x})s_{\lambda'}(\mathbf{v}),\\[5pt]
         \cdots B(x_2)B(x_1)=\sum_{\lambda }  s_{\lambda}(\mathbf{x})s_{\lambda}(\mathbf{v}).
    \end{align*}
Here, $s_{\lambda}(\mathbf{x})$ is the ordinary Schur function in commutative variables, and for the noncommutative analogue $s_{\lambda}(\mathbf{v})$,  see \cite{Fomin} for the precise definition.
\end{itemize}

\end{theorem}

\begin{theorem}\label{BGU-09}
The operators   $ \tilde{u}_i  $ satisfy the  relations in Theorem  \ref{F}.
\end{theorem}

Before we present a proof of Theorem \ref{BGU-09}, we need  two lemmas. 

    \begin{lemma} \label{u}
 The operators   $u_i  $ satisfy the   conditions in Theorem  \ref{F}. 
\end{lemma}

\begin{myproof}
When $\alpha=1$, it was verified in   \cite{Lam} that $u_i|_{\alpha=1}$ satisfy the   relations in Theorem \ref{F}. There is no obstruction to applying the arguments  to the $\alpha$ case. 
\end{myproof}
 
The second lemma lists several  relations satisfied by $u_i$ and $d_{\mu,i}$. 

\begin{lemma}\label{cl}
Let us use $d_i$ to represent  $d_{\mu,i}$ for simplicity. Then \begin{itemize}
    \item[(1)] $u_id_iu_i=\alpha\beta u_{i}^2+\beta u_i$;

    \item[(2)]  $u_id_j=d_ju_i, \forall\ i\ne j$;

     \item[(3)] $u_{i+1}d_id_{i+1}= u_{i+1}d_{i+1}d_i$;

     \item[(4)] $u_iu_{i+1}d_iu_i= \alpha\beta u_iu_{i+1}u_i+\beta u_iu_{i+1}$;

     \item[(5)] $u_{i+1}d_{i+1}u_iu_{i+1}= \alpha\beta u_{i+1}u_iu_{i+1}+\beta u_iu_{i+1}$.
\end{itemize}  
\end{lemma}

\begin{myproof}
(1) It is trivial when $u_i\cdot\lambda=0$ (namely, $i>1$ and $\lambda'_{i-1}=\lambda'_i$). Now we assume that  $u_i\cdot\lambda \ne 0$.  When $i=1$,   it is easy to see that  \[(d_iu_i-\alpha \beta u_i)\cdot\lambda=\beta  \lambda,\] 
which implies that  $u_i(d_iu_i-\alpha \beta u_i)=\beta u_i$. 
When $i\ne 1$ (in this case, $\lambda'_{i-1}>\lambda'_i$),  we can check that 
\[(d_iu_i-\alpha \beta u_i)\cdot\lambda=\beta \lambda-\alpha^k \beta  \nu,\]
where $k=\lambda'_{i-1}-\lambda'_i$, and  $\nu$ is the partition obtained from $\lambda$ by adding $k$ boxes in the $i$-th column.  Noticing that $u_i\cdot \nu=0$, we still have  
$u_i(d_iu_i-\alpha \beta u_i)=\beta u_i$. This verifies that relation in (1). 

(2)  This is  obvious by the definitions of $u_i$ and $d_j$.
       
 (3)  Notice that     $  d_id_{i+1}-d_{i+1}d_i=0$ unless $\lambda$ has the same number of boxes in the $i$-th and $(i+1)$-th columns, but $u_{i+1}$ will turn this case to 0. This leads us to $u_{i+1}(d_id_{i+1}-d_{i+1}d_i)=0$, as required. 

 (4) This case can be checked   using similar arguments as in (1). 

(5)   We only need to consider the case when both $u_i\cdot\lambda \ne 0 $ and $ u_{i+1}\cdot\lambda \ne 0$. By the analysis in the proof of (1), we have 
\[(d_{i+1}u_{i+1}-\alpha \beta u_{i+1})\cdot\lambda=\beta \lambda-\alpha^k \beta \nu,\] 
where  $k=\lambda'_{i-1}-\lambda'_{i}$, and  $\nu$ is the partition obtained from $\lambda$ by adding $k$ boxes in the $(i+1)$-th column.
Moreover, it is not hard to check that 
\[u_{i+1}u_i(\beta \lambda-\alpha^k \beta \nu)=\beta u_iu_{i+1}\cdot\lambda.\]
Collecting  the above, we deduce that 
\[
u_{i+1}u_i(d_{i+1}u_{i+1}-\alpha\beta u_{i+1})=\beta u_iu_{i+1},    
\]
which gives 
\[\alpha\beta u_{i+1}u_iu_{i+1}+\beta u_iu_{i+1}=u_{i+1}u_id_{i+1}u_{i+1}=u_{i+1}d_{i+1}u_iu_{i+1},\]
where the second equality used the relation in (2). 
    \end{myproof}
    
We can now provide a proof of 
Theorem \ref{BGU-09}. 

\noindent
{\it Proof of Theorem \ref{BGU-09}.} 
Notice that  $\tilde{u}_i $ and $\tilde{u}_k$ commute when $ | i-j  | > 1$. So the relation in (i) and (ii) are true. It remains to verify the relation in (iii). This   is trivial    when $j-i>1$ .
We next focus on  the case  when $j=i+1$. 

First, we   deduce that 
\begin{align*}
            \tilde{u}_i \tilde{u}_j-\tilde{u}_j\tilde{u}_i&=(u_iu_j-u_ju_i)(1+d_i+d_j)+u_iu_jd_id_j-u_ju_id_jd_i\\
            &=(u_iu_j-u_ju_i)(1+d_i+d_j)+u_iu_jd_jd_i-u_ju_id_jd_i\\
            &=(u_iu_j-u_ju_i)(1+d_i+d_j+d_jd_i),  
\end{align*}
where   the second equality used (3) in Lemma \ref{cl}.
So we have 
\begin{align}
  \tilde{u}_j    (\tilde{u}_i \tilde{u}_j-\tilde{u}_j\tilde{u}_i)&=(u_j+u_jd_j)(u_iu_j-u_ju_i)(1+d_i+d_j+d_jd_i)\nonumber\\
          &=u_j(u_iu_j-u_ju_i)(1+d_i+d_j+d_jd_i)\nonumber\\
         &\ \ \ \ +u_jd_ju_iu_j(1+d_i+d_j+d_jd_i) \nonumber\\
&\ \ \ \ -u_jd_ju_ju_i(1+d_i+d_j+d_jd_i).\label{FY-98}   
\end{align}
Invoking  (5) in Lemma \ref{cl}, \eqref{FY-98} can be rewritten as 
 \begin{align}
  \tilde{u}_j    (\tilde{u}_i \tilde{u}_j-\tilde{u}_j\tilde{u}_i)&=u_j(u_iu_j-u_ju_i)(1+d_i+d_j+d_jd_i)\nonumber\\
 &\ \ \ \ +(\alpha \beta u_ju_iu_j+\beta u_iu_j)(1+d_i+d_j+d_jd_i)\nonumber\\
         &\ \ \ \ -u_jd_ju_ju_i(1+d_i+d_j+d_jd_i),\nonumber
\end{align}
which, along with (1) in Lemma \ref{cl}, further gives 
 \begin{align}
 \tilde{u}_j    (\tilde{u}_i \tilde{u}_j-\tilde{u}_j\tilde{u}_i) &=u_j(u_iu_j-u_ju_i)(1+d_i+d_j+d_jd_i)\nonumber\\
  &\ \ \ \ +(\alpha \beta u_ju_iu_j+\beta u_iu_j)(1+d_i+d_j+d_jd_i)\nonumber\\
         &\ \ \ \ -(\alpha \beta u_ju_ju_i+\beta u_ju_i)(1+d_i+d_j+d_jd_i). \label{ft-34}
\end{align}
We can   write \eqref{ft-34} as
 \begin{align}
  \tilde{u}_j    (\tilde{u}_i \tilde{u}_j-\tilde{u}_j\tilde{u}_i)&=u_j(u_iu_j-u_ju_i)(1+d_i+d_j+d_jd_i)\nonumber\\
   & \ \ \ \ +
          \alpha \beta u_j(u_iu_j-u_ju_i)(1+d_i+d_j+d_jd_i)\nonumber\\
         &\ \ \ \ + \beta(u_iu_j-u_ju_i)(1+d_i+d_j+d_jd_i).\label{qsxc-2}
\end{align}

We now compute $(\tilde{u}_i \tilde{u}_j-\tilde{u}_j\tilde{u}_i)  \tilde{u}_i$. 
By direct calculation, we have 
\begin{align}
     &(\tilde{u}_i \tilde{u}_j-\tilde{u}_j\tilde{u}_i)  \tilde{u}_i =(u_iu_j-u_ju_i)(1+d_j)(1+d_i)u_i(1+d_i)\nonumber\\ 
      &=u_iu_j(1+d_j)(1+d_i)u_i(1+d_i)-u_ju_i(1+d_j)(1+d_i)u_i(1+d_i).\label{rei-09}
        \end{align}
We  compute the two terms appearing in  \eqref{rei-09} separately. Let us first consider the first term. By  (3) in Lemma \ref{cl}, we know that  for $j=i+1$,
\[u_j(1+d_j)(1+d_i)=u_j(1+d_i)(1+d_j).\] 
So the first term becomes 
\begin{equation}\label{vgu-23}
   u_iu_j(1+d_j)(1+d_i)u_i(1+d_i)=u_iu_j(1+d_i)(1+d_j)u_i(1+d_i).  
\end{equation}
We proceed to apply  (2)  and (4) in Lemma \ref{cl} to \eqref{vgu-23}, and  deduce that 
        \begin{align}
            &u_iu_j(1+d_j)(1+d_i)u_i(1+d_i)=u_iu_j(1+d_i)u_i(1+d_j)(1+d_i)\nonumber\\
           &=u_iu_ju_i(1+d_j)(1+d_i) +(\alpha \beta u_iu_ju_i+\beta u_iu_j)(1+d_j)(1+d_i).\label{YY-64}
        \end{align}
We now consider the second term in \eqref{rei-09}.  Using (1) and (2) in Lemma \ref{cl}, it follows that  
\begin{align}
           &u_ju_i(1+d_j)(1+d_i)u_i(1+d_i)=u_j(1+d_j)u_i(1+d_i)u_i(1+d_i)\nonumber\\
            &=u_j(1+d_j)u_iu_i(1+d_i)+u_j(1+d_j)(\alpha \beta u_iu_i+u_i)(1+d_i)\nonumber\\
            &=u_ju_iu_i(1+d_j)(1+d_i)+(\alpha \beta u_ju_iu_i+\beta u_ju_i)(1+d_j)(1+d_i).\label{YY-646}
\end{align}
Plugging   \eqref{YY-64} and \eqref{YY-646} into \eqref{rei-09}, we arrive at
        \begin{align}
            (\tilde{u}_i \tilde{u}_j-\tilde{u}_j\tilde{u}_i)  \tilde{u}_i
          &=(u_iu_j-u_ju_i)u_i(1+d_i+d_j+d_jd_i)\nonumber\\
          &\ \ \ \ +\alpha \beta (u_iu_j-u_ju_i)u_i(1+d_i+d_j+d_jd_i)\nonumber\\
         &\ \ \ \ +\beta(u_iu_j-u_ju_i)(1+d_i+d_j+d_jd_i)\nonumber\\
          &=u_j(u_iu_j-u_ju_i)(1+d_i+d_j+d_jd_i)\nonumber\\
          &\ \ \ \ +\alpha \beta u_j(u_iu_j-u_ju_i)(1+d_i+d_j+d_jd_i)\nonumber\\
         &\ \ \ \ +\beta(u_iu_j-u_ju_i)(1+d_i+d_j+d_jd_i),\label{uy-87}
        \end{align}
where, in  the second equality, we have used  Lemma \ref{u}.

Comparing \eqref{qsxc-2} with \eqref{uy-87} concludes  the relation 
\[\tilde{u}_j    (\tilde{u}_i \tilde{u}_j-\tilde{u}_j\tilde{u}_i)=  (\tilde{u}_i \tilde{u}_j-\tilde{u}_j\tilde{u}_i)\tilde{u}_i,\]
as required. This  completes the proof.
\qed

\begin{remark}
 Let 
\[
e_k(\mathbf{\tilde{u}})=e_k(\tilde{u}_1,\tilde{u}_2,\ldots)=\sum_{a_1>a_2>\cdots >a_k}^{} \tilde{u}_{a_1}\tilde{u}_{a_2}\cdots \tilde{u}_{a_k}. 
\] 
In light of  Theorems  \ref{F} and   \ref{BGU-09}, we know that $e_k(\mathbf{\tilde{u}})$ commute. This implies that $\tilde{A}(x)\tilde{A}(y)=\tilde{A}(y)\tilde{A}(x).$
Along with Proposition \ref{YTRd-9}, we get an alternative  proof of the symmetry   of     ${G}_{\lm}(\mathbf{x} ;\alpha;\beta)$. However, this approach  seems frustrating  in the general setting for $\mathbf{t}$ and $\mathbf{w}$. 

\end{remark}

Let 
\[
 \tilde{B}(x)=\frac{1}{1-x\tilde{u}_1}\frac{1}{1-x\tilde{u}_2}\cdots.  
\]
By Proposition \ref{YTRd-9}, we see that 
\begin{align}
  {J}_{\lm}(\mathbf{x} ;\alpha;\beta) =&  \omega \left({G}_{\lm}(\mathbf{x} ;\alpha;\beta)\right)=\omega \left \langle   \cdots \tilde{A}(x_2)\tilde{A}(x_1) \cdot \mu,\lambda \right \rangle \nonumber\\[5pt]
     =&\omega \left \langle   \left(\sum_{\lambda }  s_{\lambda}(\mathbf{x})s_{\lambda'}(\tilde{\mathbf{u}})\right)\cdot\mu,\lambda \right \rangle \nonumber=\left \langle   \left(\sum_{\lambda }  \omega (s_{\lambda}(\mathbf{x}))s_{\lambda'}(\tilde{\mathbf{u}})\right)\cdot\mu,\lambda \right \rangle\nonumber\\[5pt]
       =&\left \langle   \left(\sum_{\lambda }  s_{\lambda'}(\mathbf{x})s_{\lambda'}(\tilde{\mathbf{u}})\right)\cdot\mu,\lambda \right \rangle.\label{iue-5}
\end{align}
Combining Theorems  \ref{F} and  \ref{BGU-09}, the equality in \eqref{iue-5} becomes 
 \begin{align}\label{Uhfgh-9}      {J}_{\lm}(\mathbf{x} ;\alpha;\beta)  =\left \langle  \cdots \tilde{B}(x_2)\tilde{B}(x_1)\cdot\mu,\lambda \right \rangle.
\end{align}
Using 
\eqref{Uhfgh-9}, we can derive  a combinatorial  formula for  ${J}_{\lm}(\mathbf{x} ;\alpha;\beta)$ in terms of {\it marked multiset-valued  tableaux}.
\begin{definition}
A
marked multiset-valued  tableau
of shape $\lambda/\mu$ is a filling of the boxes of $\lambda/\mu$ with nonempty finite multisets of positive integers  such that 
\begin{itemize}
    \item The multisets are strictly increasing along each row, and weakly increasing down each column;

    \item In each box, permute  the numbers in (weakly) increasing order from left to right.  Then the leftmost  number, say   $i$, in the box can be marked with a bar if  there exists  an $i$ located at a higher position in the same column.
\end{itemize}    
\end{definition}
For a marked multiset-valued  tableau $T$,  define 
\begin{itemize}
    \item  {\it unmarkd content} $\mathrm{umcont}(T)=(n_1,n_2,\ldots)$:  $n_i$ is the number of appearances of unmarked  $i$ in $T$;

    \item {\it marked quantity vector} $\mathrm{mark}(T)=(m_1,m_2,\ldots)$: $m_i$ is the number of marked entries in the $i$-row of $T$.
\end{itemize}
Moreover, the {\it excess vector} $\mathrm{ex}(T)$ is defined in the same way as a SVRPP. 
For example, 
Figure \ref{lat-09} depicts a marked multiset-valued  tableau with
\[
  \mathrm{umcont}(T)=(3,2,2,2,1),\ \ \ \   \mathrm{mark}(T)=(0,2,1),\ \ \ \ \mathrm{ex}(T)=(3,1,1).
\]
\begin{figure}[h t]
\begin{center}
\begin{tabular}{c}
\begin{ytableau}
111& 2 & 3& 45\\
\overline{1} & 2 &\overline{3}3 \\
\overline{1}4  \\

\end{ytableau}  
\end{tabular}
\end{center}
\caption{A marked multiset-valued tableau.}\label{lat-09}
\end{figure}

\begin{remark}
In the case when $\mathrm{mark}(T)=0$ (resp., $\mathrm{ex}(T)=0$), a marked multiset-valued tableau $T$  becomes a weak set-valued tableau (resp., a valued-set tableau), as defined by Lam and 
Pylyavskyy \cite{Lam}. 
\end{remark}

Let $\mathrm{MMSVT}(\lambda/\mu)$ be the set of   
marked multiset-valued  tableaux of shape $\lambda/\mu$.

\begin{theorem}[=Theorem \ref{THM-33}]\label{ksi-01}
 For any skew shape $\lambda/\mu$, we have 
 \[
 {J}_{\lm}(\mathbf{x} ;\alpha;\beta)=\sum_{T\in \mathrm{MMSVT}(\lambda/\mu)  }   \alpha^{|\mathrm{mark}(T)|}\beta^{|\mathrm{ex}(T)|} \mathbf{x}^{\mathrm{umcont}(T)}.
 \]
\end{theorem}

\begin{myproof}
Recall that 
\[
 \tilde{B}(x_j)=\frac{1}{1-x_j\tilde{u}_1}\frac{1}{1-x_j\tilde{u}_2}\cdots.  
\]
So, by \eqref{Uhfgh-9}, we need to consider the effect of $x_j^k\tilde{u}_i^k=x_j^k(u_i+u_id_{\mu,i})^k$, which acts on the $i$-th column of a partition. The whole procedure includes $k$ steps,  and in each step, we apply either $u_i$ or $u_id_{\mu,i}$.

If we apply $u_i$, then we may add any possible    vertical strip to the $i$-th column. In this case, we fill each box in the added strip with $j$, and then  mark all the  filled $j$'s, except for the one at the  highest position, with a bar. If the added strip has $a$ boxes, then this step will contribute a factor $\alpha^{a-1}x_j$. Note that $a-1$ is just the number of  marked entries. 

If we apply $u_id_{\mu,i}$, then we first ignore the bottom box (if removable) in the $i$-th column, and then add any possible   vertical strip to the $i$-th column (so the removable box is actually not removed). 
In this case, we fill each box in the added  strip with $j$, and then  mark all the filled   $j$'s, except for the one at the highest position, with a bar. Note that this step will produce an extra (possibly repeated) number in
the removable box. Still, if the added strip has $a$ boxes, then this step will contribute a factor $\alpha^{a-1}\beta x_j$.

Let $T$ be  a filling of $\lambda/\mu$ generated in the procedure $\left \langle  \cdots \tilde{B}(x_2)\tilde{B}(x_1)\cdot\mu,\lambda \right \rangle$. By the above analysis, each box may be filled with a nonempty finite multiset, and the multisets are weakly increasing in each column. Notice that  in $ \tilde{B}(x_j)$, the actions  on columns are   from right to left. So the multisets in each row will be strictly increasing. 
Moreover, by the above construction, it is easy to see that  $T$ satisfies the restrictions imposed upon the barred entries as in the definition of a marked multiset-valued tableau. So $T$ belongs to  $\mathrm{MMSVT}(\lambda/\mu)$. 

In the reverse direction, it is not hard to check that each $T\in \mathrm{MMSVT}(\lambda/\mu)$ can be (uniquely) produced via the above procedure.  Furthermore, each such produced  $T$ has an associated weight $\alpha^{|\mathrm{mark}(T)|}\beta^{|\mathrm{ex}(T)|} \mathbf{x}^{\mathrm{umcont}(T)}$.
This finishes the proof.
\end{myproof}
 
Setting $\alpha=0$ and $\beta=-1$ (resp., $\alpha=1$ and $\beta=0$),   Theorem  \ref{ksi-01} recovers the image of a stable Grothendieck polynomial (or, a dual stable Grothendieck polynomial) as given  in \cite{Lam}. 


\section{Concluding remarks}\label{Last-22}

In this section, we discuss some problems and 
conjectures concerning hybrid Grothendieck polynomials. We shall mainly concentrate on  
\[
    {G}_{\lm}(\mathbf{x} ;\alpha;\beta)={G}_{\lm}(\mathbf{x} ;\mathbf{t};\mathbf{w})|_{t_i=\alpha, w_i=\beta}.
\]
Keep in mind that ${G}_{\lm}(\mathbf{x} ;\mathbf{0};-\mathbf{1})$  and  ${G}_{\lm}(\mathbf{x} ;\mathbf{1};\mathbf{0})$ are equal to the   stable and dual stable  Grothendieck polynomials, respectively. 
 
\subsection{The basis problem}
Consider the ring 
 \[\Lambda=\bigoplus_{\lambda}\mathbb{Q}(\alpha,\beta) s_{\lambda}(\mathbf{x})\]  of symmetric functions over the field of rational functions in $\alpha$ and $\beta$. Denote by $\hat \Lambda$  the completion of $\Lambda$ which may contain  infinite linear combinations of Schur functions. In other words, $\hat \Lambda$ consists of formal power series possibly with unbounded degree.  
 
Both  stable and dual stable Grothendieck polynomials form a basis of $\hat \Lambda$  since the corresponding transition matrices, with respect to the basis of Schur function, are triangular \cite{Lam}. We conjecture that 
  hybrid Grothendieck polynomials also form a basis.

 \begin{conjecture}\label{conj-3-2}
The family of  hybrid Grothendieck polynomial ${G}_{\lambda}(\mathbf{x};\alpha;\beta)$ constitute  a basis of $\hat \Lambda$.   
\end{conjecture}

There have been  Pieri type formulas for stable or dual stable Grothendieck polynomials \cite{Iwao1,Iwao2,Lenart,Yeliussizov 1,Yeliussizov 2}, that is, a formula for the expansion of a stable or dual stable Grothendieck polynomial by multiplying a Schur function (or a stable/dual stable Grothendieck polynomial) indexed by a one-row or one-column partition.

 \begin{problem}
Assuming Conjecture  \ref{conj-3-2}, do there exist  Pieri type formulas for   ${G}_{\lambda}(\mathbf{x} ;\alpha;\beta)$?   
\end{problem}

\subsection{Determinantal   formulas}

Let us recall the following two kinds of determinantal formulas for Schur polynomials. The first one can be  described via Cauchy's bialternant formula, which is also called Weyl's formula,
\[
s_{\lambda}(\mathbf{x}_n)=\frac{\det\left(x_i^{\lambda_j+n-j}\right)_{i,j=1}^n}{\prod_{1\leq i<j\leq n}(x_i-x_j)}.
\]
The second one  is the Jacobi--Trudi type identity which applies to   skew shapes:
\[
s_{\lambda/\mu}(\mathbf{x}_n)=\det \left(h_{\lambda_i-\mu_j-i+j}(\mathbf{x}_n)\right)_{i,j=1}^n=\det \left(e_{\lambda_i'-\mu_j'-i+j}(\mathbf{x}_n)\right)_{i,j=1}^n,
\]
where $e_k(\mathbf{x}_n)$ and $h_k(\mathbf{x}_n)$ are respectively the elementary and complete symmetric polynomials: 
\[
 e_k(\mathbf{x}_n)=\sum_{1\leq i_1<\cdots<i_k\leq n}x_{i_1}\cdots x_{i_k}, \ \ \ \     h_k(\mathbf{x}_n)=\sum_{1\leq i_1\leq \cdots\leq i_k\leq n}x_{i_1}\cdots x_{i_k}.
\]


The Weyl type formulas for stable Grothendieck polynomials {or dual stable Grothendieck polynomials} have been given  for example  in \cite{Amanov,Ikeda-1,Ikeda-2,McNamara}. 
On the other hand, the Jacobi--Trudi type formulas for stable or dual stable Grothendieck polynomials (more generally for canonical or refined canonical versions) have been extensively explored, see for example  \cite{Amanov,Hwang,Kim1,Lascoux2,Matsumura1,Matsumura2,Yeliussizov 1}.

\begin{problem}
Does the hybrid Grothendieck polynomial ${G}_{\lm}(\mathbf{x}_n ;\alpha;\beta)$ admit  Weyl type or Jacobi--Trudi type formuals?   
\end{problem}

\subsection{Expansions in  other basis}

Other than using the Schur basis,  we may consider the expansion of a hybrid Grothendieck polynomial in the basis of stable or dual stable Grothendieck polynomials.  
Kundu \cite{Kundu} obtained the expansion of the refined stable Grothendieck polynomial $G_{\lambda/\mu}(\mathbf{x};\mathbf{0};\mathbf{w})$ in the basis $\{G_{\nu}(\mathbf{x};\mathbf{0};\mathbf{-1})\}$ of stable  Grothendieck polynomials as well as  in the basis $\{G_{\nu}(\mathbf{x};\mathbf{1};\mathbf{0})\}$ of dual stable  Grothendieck polynomials. Similar expansions for refined dual stable Grothendieck polynomials $G_{\lambda/\mu}(\mathbf{x};\mathbf{t};\mathbf{0})$ were also given in \cite{Kundu}. The proofs used the `tableaux-Schur
expansion' method due to  Bandlow and  Morse \cite{Bandlow}. For related  expansions for  canonical Grothendieck polynomials, see  the work of Pan,   Pappe,   Poh and   Schilling \cite{Pan}.

\begin{problem}
Find the expansion of a   hybrid Grothendieck polynomial   in the basis of  stable or dual stable Grothendieck polynomials.   
\end{problem}

It should   be pointed  out that the square root crystal approach  in \cite{MTY} also provides a way to decompose a symmetric function in terms of stable Grothendieck polynomials. We do not know if there exists such a crystal on set-valued reverse plane partitions. 
 
One can also define a flagged analogue  for hybrid Grothendieck polynomials. 
A flag $\phi=(\phi_1\leq \phi_2\leq \cdots)$ is a  sequence of weakly increasing positive integers. A  SVRPP flagged by $\phi$ is a SVRPP such that each number in row $i$ cannot exceed $\phi_i$.
The flagged hybrid Grothendieck polynomial ${G}_{\lm}(\mathbf{x}_\phi ;\alpha;\beta)$
is defined as the generating polynomial  over SVRPP's  flagged 
 by $\phi$.
Note that   ${G}_{\lm}(\mathbf{x}_\phi ;\alpha;\beta)$ is not necessarily symmetric. 

By building  the   Demazure cystal structures, flagged stable or flagged dual stable Grothendieck polynomials can be decomposed in terms of key polynomials \cite{Kundu12,Kundu}. 

\begin{problem}
Does there exist a Demazure crystal structure on flagged set-valued reverse plane partitions? 
\end{problem}

\subsection{Integrable lattice model}

Integeral lattice models from statistical mechanics  have been a powerful tool in the study of Schubert calculus as well as the properties of important families of symmetric/asymmetric functions in combinatorics.   See for example \cite{ABW,Buciumas,FGX,Wheeler} and the references theirin.

Different  vertex models for stable or dual stable  Grothendieck
polynomials have appeared for example in the work of   Gunna and  Zinn-Justin
  \cite{Gunna} and Motegi and Scrimshaw \cite{Motegi}.
 
\begin{problem}
Is there an integral lattice model for the hybrid Grothendieck polynomial ${G}_{\lm}(\mathbf{x}_n ;\alpha;\beta)$?   
\end{problem}

\subsection{Algorithm for the decomposition into Schur functions}

Consider the Schur expansion of the stable Grothendieck polynomial
${G}_{\lambda}(\mathbf{x} ;\mathbf{0};-\mathbf{1})$  or the dual stable Grothendieck polynomial ${G}_{\lambda}(\mathbf{x} ;\mathbf{1};\mathbf{0})$ of straight shape. 
As mentioned in Introduction, the Schur expansion of 
${G}_{\lambda}(\mathbf{x} ;\mathbf{0};-\mathbf{1})$ was first given by  Lenart \cite{Lenart} where the coefficients were interpreted  in terms of strict tableaux. The proof given in \cite{Lenart} is not purely combinatorial. 
Buch \cite{Buch} found a combinatorial proof by developing  the uncrowding algorithm, which maps a set-valued tableau into a pair of    semistandard Young tableau and   strict tableau.  See  \cite{Pan} for an extension of the uncrowding algorithm to hook-valued tableaux.  

On the other hand, 
Lam and Pylyavskyy \cite{Lam} proved the Schur expansion of   the dual stable Grothendieck polynomial ${G}_{\lambda}(\mathbf{x} ;\mathbf{1};\mathbf{0})$ where the coefficients were described by   elegant fillings. The proof in \cite{Lam} is based on an algorithm mapping a reverse plane partition into a pair of semistandard Young tableau and  elegant filling. 

We  beg for a combinatorial proof for the Schur expansion of a hybrid  Grothendieck polynomial.

\begin{problem}
Find a combinatorial algorithm to prove the Schur expansion of the hybrid  Grothendieck polynomial ${G}_{\lambda}(\mathbf{x};\alpha;\beta)$ of straight shape.  
\end{problem}

\subsection{Coincidences among hybrid Grothendieck polynomials }

The coincidence problem about Schur functions is to determine when two
different skew shapes give rise to the same skew Schur function.   For progress on this problem, see for example \cite{Billera,JL,McNamara-2,Reiner}. The   coincidences  among stable or dual stable Grothendieck polynomials  have been considered  by  Alwaise et al.  \cite{Alwaise}.

\begin{problem}
Characterize the  skew shapes $\lambda/\mu$ and $\nu/\rho$ that produce the same     hybrid Grothendieck   polynomial, namely,   $ {G}_{\lm}(\mathbf{x} ;\alpha;\beta)={G}_{\nu/\rho}(\mathbf{x} ;\alpha;\beta)$.   
\end{problem}

Let $\delta_n=(n,n-1,\ldots,1)$ denote the staircase partition. The Stembridge equality   states that  \[s_{\delta_n/\rho}(\mathbf{x})=s_{\delta_n/\rho'}(\mathbf{x})\]
holds for any $\rho\subseteq \delta_n$, where $\rho'$ is the conjugate partiton of $\rho$.
Analogues of the Stembridge equality for     stable and   dual stable Grothendieck polynomials   have been given  by Abney-McPeek,  An and Ng \cite{Abney-McPeek}, that is,  
$$G_{\delta_n/\rho}(\mathbf{x};\mathbf{0};-\mathbf{1})=G_{\delta_n/\rho'}(\mathbf{x};\mathbf{0};-\mathbf{1}),\quad\quad
G_{\delta_n/\rho}(\mathbf{x};\mathbf{1};\mathbf{0})=G_{\delta_n/\rho'}(\mathbf{x};\mathbf{1};\mathbf{0}).$$
Note that the converse   is also true \cite{Abney-McPeek}:  if  
\begin{equation}\label{njs-9-31}
G_{\nu/\rho}(\mathbf{x};\mathbf{0};-\mathbf{1})=G_{\nu/\rho'}(\mathbf{x};\mathbf{0};-\mathbf{1}),\ \ \ \text{or}\ \ \ 
G_{\nu/\rho}(\mathbf{x};\mathbf{1};\mathbf{0})=G_{\nu/\rho'}(\mathbf{x};\mathbf{1};\mathbf{0})     
 \end{equation} 
 holds for all $\rho\subseteq \nu$, then $\nu$ must be a  staircase partition. We conjecture the following analogue for hybrid Grothendieck polynomials. 

\begin{conjecture}
 For  any $\rho\subseteq \delta_n$, we have
 \[G_{\delta_n/\rho'}(\mathbf{x};\alpha;\beta)=G_{\delta_n/\rho'}(\mathbf{x};\alpha;\beta).\]
 \end{conjecture}

Of course, in view of \eqref{njs-9-31}, if 
 \[
 G_{\nu/\rho}(\mathbf{x};\alpha;\beta)=G_{\nu/\rho'}(\mathbf{x};\alpha;\beta) \]
 holds for all $\rho\subseteq \nu$, then $\nu$ must be a staircase partition.

\clearpage
\begin{appendices}

\section{The crystal graph of $\mathrm{SVRPP}^{n}(\lambda/\mu)$ with $\l=(2,2)$, $\m=(1)$ and $n=3$}\label{hgi-6-9}

\begin{figure}[H]
\[
\ytableausetup{boxsize=2em}
\scalebox{0.7}{
\begin{tikzpicture}[>=latex]
\node (A) {\ytableaushort{\none 1,11}};
\node[below = 1 of A] (B) {\ytableaushort{\none 2,12}};
\node[below left = 1 and 0.5 of B] (C) {\ytableaushort{\none 2,22}};
\node[below right = 1 and 0.5 of B] (D) {\ytableaushort{\none 3,13}};
\node[below right = 1 and 0.5 of C] (E) {\ytableaushort{\none 3,23}};
\node[below = 1 of E] (F) {\ytableaushort{\none 3,33}};
\path (A) edge[pil,color=blue] node[right,black]{$1$} (B);
\path (B) edge[pil,color=blue] node[above left,black]{$1$} (C);
\path (B) edge[pil,color=red] node[above right,black]{$2$} (D);
\path (C) edge[pil,color=red] node[below left,black]{$2$} (E);
\path (D) edge[pil,color=blue] node[below right,black]{$1$} (E);
\path (E) edge[pil,color=red] node[right,black]{$2$} (F);
\end{tikzpicture}
\hspace{6mm}

\begin{tikzpicture}[>=latex]
\node (A) {\ytableaushort{\none 1,1{12}}};
\node[below left = 1 and 0.5 of A] (B) {\ytableaushort{\none 2,{12}2}};
\node[below right = 1 and 0.5 of A] (C) 
{\ytableaushort{\none 1,1{13}}};
\node[below = 1 of B] (D) {\ytableaushort{\none 3,{12}3}};
\node[below = 1 of C] (E) {\ytableaushort{\none 2,1{23}}};
\node[below = 1 of D] (F) {\ytableaushort{\none 3,{13}3}};
\node[below = 1 of E] (G) {\ytableaushort{\none 2,2{23}}};
\node[below right = 1 and 0.5 of F] (H) {\ytableaushort{\none 3,{23}3}};

\path (A) edge[pil,color=blue] node[above left,black]{$1$} (B);
\path (A) edge[pil,color=red] node[above right,black]{$2$} (C);
\path (B) edge[pil,color=red] node[right,black]{$2$} (D);
\path (C) edge[pil,color=blue] node[right,black]{$1$} (E);
\path (D) edge[pil,color=red] node[right,black]{$2$} (F);
\path (E) edge[pil,color=blue] node[right,black]{$1$} (G);
\path (F) edge[pil,color=blue] node[below left,black]{$1$} (H);
\path (G) edge[pil,color=red] node[below right,black]{$2$} (H);
\end{tikzpicture}
\hspace{6mm}
\begin{tikzpicture}[>=latex]
\node (A) {\ytableaushort{\none 1,1{2}}};
\node[below left = 1 and 0.5 of A] (B) {\ytableaushort{\none 1,{2}2}};
\node[below right = 1 and 0.5 of A] (C) 
{\ytableaushort{\none 1,1{3}}};
\node[below = 1 of B] (D) {\ytableaushort{\none 1,{2}3}};
\node[below = 1 of C] (E) {\ytableaushort{\none 2,1{3}}};
\node[below = 1 of D] (F) {\ytableaushort{\none 1,{3}3}};
\node[below = 1 of E] (G) {\ytableaushort{\none 2,2{3}}};
\node[below right = 1 and 0.5 of F] (H) {\ytableaushort{\none 2,33}};

\path (A) edge[pil,color=blue] node[above left,black]{$1$} (B);
\path (A) edge[pil,color=red] node[above right,black]{$2$} (C);
\path (B) edge[pil,color=red] node[right,black]{$2$} (D);
\path (C) edge[pil,color=blue] node[right,black]{$1$} (E);
\path (D) edge[pil,color=red] node[right,black]{$2$} (F);
\path (E) edge[pil,color=blue] node[right,black]{$1$} (G);
\path (F) edge[pil,color=blue] node[below left,black]{$1$} (H);
\path (G) edge[pil,color=red] node[below right,black]{$2$} (H);
\end{tikzpicture}

}
\]
\end{figure}

\begin{figure}[H]
\[
\ytableausetup{boxsize=2em}
\scalebox{0.7}{
\begin{tikzpicture}[>=latex]
\node (A) {\ytableaushort{\none 1,{12}2}};
\node[below = 1 of A] (B) {\ytableaushort{\none 1,{12}3}};
\node[below left = 1 and 0.5 of B] (C) {\ytableaushort{\none2,{12}3}};
\node[below right = 1 and 0.5 of B] (D) {\ytableaushort{\none1,{13}3}};
\node[below right = 1 and 0.5 of C] (E) {\ytableaushort{\none2,{13}3}};
\node[below = 1 of E] (F) {\ytableaushort{\none2,{23}3}};
\path (A) edge[pil,color=red] node[right,black]{$2$} (B);
\path (B) edge[pil,color=blue] node[above left,black]{$1$} (C);
\path (B) edge[pil,color=red] node[above right,black]{$2$} (D);
\path (C) edge[pil,color=red] node[below left,black]{$2$} (E);
\path (D) edge[pil,color=blue] node[below right,black]{$1$} (E);
\path (E) edge[pil,color=blue] node[right,black]{$1$} (F);
\end{tikzpicture}
\hspace{6mm}
\begin{tikzpicture}[>=latex]
\node (A) {\ytableaushort{\none {12},1{2}}};
\node[below left = 1 and 0.5 of A] (B) {\ytableaushort{\none{12},2{2}}};
\node[below right = 1 and 0.5 of A] (C) 
{\ytableaushort{\none{13},1{3}}};
\node[below = 1 of B] (D) {\ytableaushort{\none {13},2{3}}};
\node[below = 1 of C] (E) {\ytableaushort{\none {23},1{3}}};
\node[below = 1 of D] (F) {\ytableaushort{\none {13},3{3}}};
\node[below = 1 of E] (G) {\ytableaushort{\none {23},2{3}}};
\node[below right = 1 and 0.5 of F] (H) 
{\ytableaushort{\none {23},3{3}}};

\path (A) edge[pil,color=blue] node[above left,black]{$1$} (B);
\path (A) edge[pil,color=red] node[above right,black]{$2$} (C);
\path (B) edge[pil,color=red] node[right,black]{$2$} (D);
\path (C) edge[pil,color=blue] node[right,black]{$1$} (E);
\path (D) edge[pil,color=red] node[right,black]{$2$} (F);
\path (E) edge[pil,color=blue] node[right,black]{$1$} (G);
\path (F) edge[pil,color=blue] node[below left,black]{$1$} (H);
\path (G) edge[pil,color=red] node[below right,black]{$2$} (H);
\end{tikzpicture}
\hspace{6mm}
\begin{tikzpicture}[>=latex]
\node (A) {\ytableaushort{\none {12},{12}2}};
\node[below = 1 of A] (B) {\ytableaushort{\none {13},{12}3}};
\node[below left = 1 and 0.5 of B] (C) 
{\ytableaushort{\none{23},{12}3}};
\node[below right = 1 and 0.5 of B] (D) 
{\ytableaushort{\none{13},{13}3}};
\node[below right = 1 and 0.5 of C] (E) 
{\ytableaushort{\none{23},{13}3}};
\node[below = 1 of E] (F) {\ytableaushort{\none{23},{23}3}};
\path (A) edge[pil,color=red] node[right,black]{$2$} (B);
\path (B) edge[pil,color=blue] node[above left,black]{$1$} (C);
\path (B) edge[pil,color=red] node[above right,black]{$2$} (D);
\path (C) edge[pil,color=red] node[below left,black]{$2$} (E);
\path (D) edge[pil,color=blue] node[below right,black]{$1$} (E);
\path (E) edge[pil,color=blue] node[right,black]{$1$} (F);
\end{tikzpicture}
}
\]
\end{figure}

\begin{figure}[H]
\[
\ytableausetup{boxsize=2em}
\scalebox{0.7}{
\begin{tikzpicture}[>=latex]
\node (A) {\ytableaushort{\none {1},1{123}}};
\node[below = 1 of A] (B) {\ytableaushort{\none {2},{12}{23}}};
\node[below = 1 of B] (C) {\ytableaushort{\none {3},{123}3}};
\path (A) edge[pil,color=blue] node[above left,black]{$1$} (B);
\path (B) edge[pil,color=red] node[right,black]{$2$} (C);
\end{tikzpicture}
\hspace{8mm}
\begin{tikzpicture}[>=latex]
\node (A) {\ytableaushort{\none {1},1{23}}};
\node[below = 1 of A] (B) {\ytableaushort{\none {1},{2}{23}}};
\node[below = 1 of B] (C) {\ytableaushort{\none {1},{23}3}};
\path (A) edge[pil,color=blue] node[above left,black]{$1$} (B);
\path (B) edge[pil,color=red] node[right,black]{$2$} (C);
\end{tikzpicture}

\hspace{8mm}
\begin{tikzpicture}[>=latex]
\node (A) {\ytableaushort{\none {12},1{23}}};
\node[below = 1 of A] (B) {\ytableaushort{\none {12},{2}{23}}};
\node[below = 1 of B] (C) {\ytableaushort{\none {13},{23}3}};
\path (A) edge[pil,color=blue] node[above left,black]{$1$} (B);
\path (B) edge[pil,color=red] node[right,black]{$2$} (C);
\end{tikzpicture}

\hspace{8mm}
\begin{tikzpicture}[>=latex]
\node (A) {\ytableaushort{\none {12},1{3}}};
\node[below = 1 of A] (B) {\ytableaushort{\none {12},{2}{3}}};
\node[below = 1 of B] (C) {\ytableaushort{\none {13},{3}3}};
\path (A) edge[pil,color=blue] node[above left,black]{$1$} (B);
\path (B) edge[pil,color=red] node[right,black]{$2$} (C);
\end{tikzpicture}

\hspace{8mm}
\begin{tikzpicture}[>=latex]
\node (A) {\ytableaushort{\none {123},1{3}}};
\node[below = 1 of A] (B) {\ytableaushort{\none {123},{2}{3}}};
\node[below = 1 of B] (C) {\ytableaushort{\none {123},{3}3}};
\path (A) edge[pil,color=blue] node[above left,black]{$1$} (B);
\path (B) edge[pil,color=red] node[right,black]{$2$} (C);
\end{tikzpicture}
}
\]
\end{figure}

\begin{figure}[H]
\[
\ytableausetup{boxsize=2em}
\scalebox{0.7}{
\begin{tikzpicture}[>=latex]
\node (A) {\ytableaushort{\none {1},{12}{23}}};
\node[below = 1 of A] (B) {\ytableaushort{\none {1},{123}{3}}};
\node[below = 1 of B] (C) {\ytableaushort{\none {2},{123}3}};
\path (A) edge[pil,color=red] node[above left,black]{$2$} (B);
\path (B) edge[pil,color=blue] node[right,black]{$1$} (C);
\end{tikzpicture}

\hspace{8mm}
\begin{tikzpicture}[>=latex]
\node (A) {\ytableaushort{\none {12},{12}{23}}};
\node[below = 1 of A] (B) {\ytableaushort{\none {13},{123}{3}}};
\node[below = 1 of B] (C) {\ytableaushort{\none {23},{123}3}};
\path (A) edge[pil,color=red] node[above left,black]{$2$} (B);
\path (B) edge[pil,color=blue] node[right,black]{$1$} (C);
\end{tikzpicture}

\hspace{8mm}
\begin{tikzpicture}[>=latex]
\node (A) {\ytableaushort{\none {12},{12}{3}}};
\node[below = 1 of A] (B) {\ytableaushort{\none {12},{13}{3}}};
\node[below = 1 of B] (C) {\ytableaushort{\none {12},{23}3}};
\path (A) edge[pil,color=red] node[above left,black]{$2$} (B);
\path (B) edge[pil,color=blue] node[right,black]{$1$} (C);
\end{tikzpicture}

\hspace{8mm}
\begin{tikzpicture}[>=latex]
\node (A) {\ytableaushort{\none {123},{12}{3}}};
\node[below = 1 of A] (B) {\ytableaushort{\none {123},{13}{3}}};
\node[below = 1 of B] (C) {\ytableaushort{\none {123},{23}3}};
\path (A) edge[pil,color=red] node[above left,black]{$2$} (B);
\path (B) edge[pil,color=blue] node[right,black]{$1$} (C);
\end{tikzpicture}

\hspace{8mm}
\begin{tikzpicture}[>=latex]
\node (A) {\ytableaushort{\none {12},{123}{3}}};
\node[below = 4 of A] (B) {\ytableaushort{\none {123},{123}{3}}};
\end{tikzpicture}
}
\]
\end{figure}

\end{appendices}

\footnotesize{

\textsc{\{Peter L. Guo, Mingyang Kang, Jiaji Liu\} Center for Combinatorics, Nankai University, LPMC, Tianjin 300071, P.R. China}

{\it
Email address: \tt lguo@nankai.edu.cn, 2120240005@mail.nankai.edu.cn,\\ 1120240006@mail.nankai.edu.cn}

\end{document}